%% file: iclr2025_conference.tex
\pdfoutput=1
\documentclass{article} 
\usepackage{iclr2025_conference,times}

\usepackage{hyperref}
\usepackage[english]{babel}
\usepackage[T1]{fontenc}
\usepackage[utf8]{inputenc}
\input{math_commands.tex}

\input{notation.tex}
\usepackage{url}

\title{Interference Among First-Price Pacing Equilibria: A Bias and Variance Analysis}

\author{  \textbf{Luofeng Liao} \quad \textbf{Christian Kroer}\\
  Columbia University \\
  \texttt{\{ll3530, ck2945\}@columbia.edu}
  \And
  \textbf{Sergei Leonenkov}\\
  Ads Online Experimentation, Meta \\
  \texttt{leonenkov@meta.com}
    \AND
  \textbf{Okke Schrijvers} \quad \textbf{Liang Shi} \quad
  \textbf{Nicolas Stier-Moses} \quad \textbf{Congshan Zhang}\\
  Central Applied Science, Meta \\
  \texttt{\{okke, liangshi, nstier, cszhang\}@meta.com} 
  }

\iclrfinalcopy 
\begin{document}

\maketitle

\begin{abstract}
A/B testing is widely used in the internet industry. For online marketplaces (such as advertising markets), standard approaches to A/B testing may lead to biased results when buyers have budget constraints, as budget consumption in one arm of the experiment impacts performance of the other arm. 
This is often addressed using a budget-split design. Yet such splitting may degrade statistical performance as budgets become too small in each arm.
We propose a \emph{parallel budget-controlled A/B testing} design where we use market segmentation to identify submarkets in the larger market, and we run parallel budget-split experiments in each submarket.
We demonstrate the effectiveness of this approach on real experiments on advertising markets at Meta.
Then, we formally study interference that derives from such experimental designs, using the first-price pacing equilibrium framework as our model of market equilibration.
We propose a debiased surrogate that eliminates the first-order bias of FPPE, and derive a plug-in estimator for the surrogate and establish its asymptotic normality. We then provide an estimation procedure for submarket parallel budget-controlled A/B tests. Finally, we present numerical examples on semi-synthetic data, confirming that the debiasing technique achieves the desired coverage properties.
\end{abstract}
\input{maintext.tex}
\bibliography{refs.bib}
\bibliographystyle{iclr2025_conference}
\appendix
\input{app_theory.tex}
\end{document}

%% file: math_commands.tex
\usepackage{amsmath,amsfonts,bm}

\def\eqref#1{equation~\ref{#1}}

\def\1{\bm{1}}

\DeclareMathAlphabet{\mathsfit}{\encodingdefault}{\sfdefault}{m}{sl}
\SetMathAlphabet{\mathsfit}{bold}{\encodingdefault}{\sfdefault}{bx}{n}

\newcommand{\E}{\mathbb{E}}

\newcommand{\R}{\mathbb{R}}

\DeclareMathOperator*{\argmax}{arg\,max}
\DeclareMathOperator*{\argmin}{arg\,min}

%% file: notation.tex
\usepackage{mathtools}   
\usepackage{amsmath}
\usepackage{amsthm}
\usepackage{multirow}
\usepackage{mathrsfs}

\usepackage{amsthm}
\usepackage[ruled]{algorithm2e} 
\usepackage{subfigure} 

\usepackage[shortlabels]{enumitem}
\setlist[itemize]{topsep=0pt,}
\usepackage{pifont}
\hypersetup{
colorlinks   = true, 
urlcolor     = blue, 
linkcolor    = blue, 
citecolor   = blue 
}

\usepackage{dsfont} 
\usepackage{xcolor} 
\usepackage{cleveref}
\usepackage{esvect}

\usepackage{tabularx}
\usepackage{booktabs}
\usepackage[labelfont=bf,format=plain,justification=raggedright,singlelinecheck=false]{caption}

\usepackage{bm}
\usepackage{calrsfs}

\newcommand{\ubar}[1]{\underline{#1}}
\newcommand{\pmr}{{pacing multiplier }}

\def \tr{\mathop{\text{tr}}\kern.2ex}

\def \sa {_\alpha}

\def\R{{\mathbb R}}
\def\Rp{{\mathbb R _+}}
\def\Rn{{\mathbb R^n}}

\def\P{{\mathbb P}}
\def\E{{\mathbb E}}
\def\B{{\mathbb B}}

\def\cH{{\boldsymbol{H}}}

\def\cM{{\mathscr{M}}}

\def\cN{{\mathcal{N}}}

\newtheorem{lemma}{Lemma}
\newtheorem{defn}{Definition}

\newlist{enumconditions}{enumerate}{1} 
\setlist[enumconditions]{label = \thelemma.\alph*}
\crefname{enumconditionsi}{Cond.}{Conditions}

\newlist{enumlmresult}{enumerate}{1} 
\setlist[enumlmresult]{label = \thelemma.\arabic*}
\crefname{enumlmresulti}{Part}{Parts}

\crefname{figure}{Fig}{Figures}
\crefname{algorithm}{Algorithm}{Algorithms}
\crefname{table}{Table}{Tables}

\newlist{enumdef}{enumerate}{1} 
\setlist[enumdef]{label = \thedefn.\arabic*}
\crefname{enumdefi}{Cond.}{Cond.}

\crefname{Assumption}{Assumption}{Assumptions}
\renewcommand*{\theAssumption}{\arabic{Assumption}}
\crefname{defn}{Def.}{Defs.}
\crefname{example}{Ex.}{Examples}

\newlist{enumthmresult}{enumerate}{1} 
\setlist[enumthmresult]{label = \thetheorem.\arabic*}
\crefname{enumthmresulti}{Part}{Parts}

\crefname{equation}{Eq.}{Eqs.}

\newlist{enumassumption}{enumerate}{1} 
\setlist[enumassumption]{label = \theAssumption.\arabic*}
\crefname{enumassumptioni}{Condition}{Conditions}

\crefname{appendix}{App.}{Apps}
\crefname{theorem}{Thm.}{Theorems}
\crefname{section}{Sec.}{Sections}

\newtheorem{theorem}{Theorem}
\renewcommand*{\thetheorem}{\arabic{theorem}}

\def \cd {{ \cdot }}
\def \deltasti {{\delta^*_i}}

\def \Ic {{I^{c}}}

\def \IF { \bm{d}}

\def \REV {{ {{REV}}}}
\def \REVst {{ {{REV}}  }^* }

\def \FPPE {{  \mathsf{{FPPE}}}}
\def \oFPPE {{  \widehat{ \textsf{{FPPE}}}}}
\def \bidgap {{  \mathsf{{bidgap}}}}

\def \CI {{  \mathbb{{CI}}}}

\def \b {{\boldsymbol  {\bm \beta}}}

\def \mutau {\bm \mu^\tau}

\def \vithetau {{ v_i(\theta^\tau) }}
\def \vithe {{ v_i(\theta) }}
\def \vitau {{ v_i^\tau }}
\def \vbar {{ \bar v }}

\def \var {{  {Var} }}
\def \cov {{  {Cov} }}

\def \sumtau {{\sum_{\tau=1}^{t}}}

\def \sumiton {{ \sum_{i=1}^n }}

\def \ptau {{ p ^ \tau }}

\def \t {\theta}
\def \thetau {{ \theta^\tau }}

\def \pst {p^*}

\def \betast {{\bm \beta}^*} 
\def \betasta {{\bm \beta}^*_\alpha}
\def \betastai {\beta^*_{\alpha , i } }
\def \betastt {\widetilde{{\bm \beta}}^*}

\def \betasti {{\beta^*_i}}

\def \betai {{\beta_i}}
\def \betak {{\beta_k}}
\def \betagam {{\bm \beta}^\gamma}
\def \betagama {{\bm \beta}^\gamma_\alpha}

\def \betagamai {\beta^\gamma_{\alpha, i}}

\def \xitau {{x_i^\tau}}

\def \xsti {{x^*_i}}

\def \good {\mathsf{good}}
\def \bad {\mathsf{bad}}
\def \hinv { ( {\bm P} \sa \cH _{\alpha} {\bm P} \sa) \pinv }
\def \hhinv { (\widehat {\bm P} \sa \widehat \cH _{\alpha} \widehat {\bm P} _ \alpha ) \pinv}
\def \inv {^{{-1}}}

\def \pinv {^{\dagger}}
\def \sq {^{2}}
\def \tp {^{\scriptscriptstyle\mathsf{T}}}
\def \st {^{*}}

\def \Diag {{ \operatorname*{ Diag}}}

\newcommand{\defeq}{=}

\newcommand*\diff{\mathop{}\!\mathrm{d}}

\def \toprob {{ \,\overset{{p}}{\to}\, }}

\def \tod {{  \,\overset{{d}}{\to}\, }}

\def \biasformula { \diff {\bm \beta}^*(\alpha)}

%% file: maintext.tex
\section{Introduction}\label{sec:intro}
Online A/B testing is widely used in the internet industry to inform decisions on new feature roll-outs. For online marketplaces (such as advertising markets), standard approaches to A/B testing may lead to biased results when buyers operate under a budget constraint, as budget consumption in one arm of the experiment impacts performance of the other arm. To counteract this interference, one can use a budget-split design where the budget constraint operates on a per-arm basis and each arm receives an equal fraction of the budget, leading to ``budget-controlled A/B testing,'' see e.g. \cite{Basse2016,Liu2021}. 

Despite clear advantages of budget-controlled A/B testing, companies are extremely constrained by the number of such experiments they can run. While it's possible to create more budget splits, this will lower the budget per group substantially, which could lead to different equilibrium outcomes and may disproportionately affect smaller buyers. Additionally, a common approach to increase experimentation throughput is to run orthogonal experiments (with their own orthogonal randomization), but this would either suffer from the same interference as the vanilla A/B test setup, or also require further budget splits.

In this paper, we propose a \emph{parallel budget-controlled A/B test} design where we use market segmentation to identify submarkets in the larger market, and we run parallel experiments on each submarket. When the overall market can be divided into several relatively isolated submarkets, budget-controlled A/B tests can be conducted in parallel within these submarkets. However, this method also presents some challenges. First, submarkets are rarely completely isolated; certain items may attract buyers from multiple submarkets, resulting in interference across 
submarkets when conducting tests in parallel. Second, submarkets differ in terms of buyer (and user) composition, which might cause the local treatment effect estimates to not be representative of the global treatment effect where all buyers are included in the market. The second challenge is relatively easy to address in practice by imposing balancing constraints in the clustering algorithm used to define submarkets, while the first challenge is more fundamentally important and requires deeper understanding. 

\begin{figure}[b]
	\centering
	\includegraphics[width=0.34\linewidth]{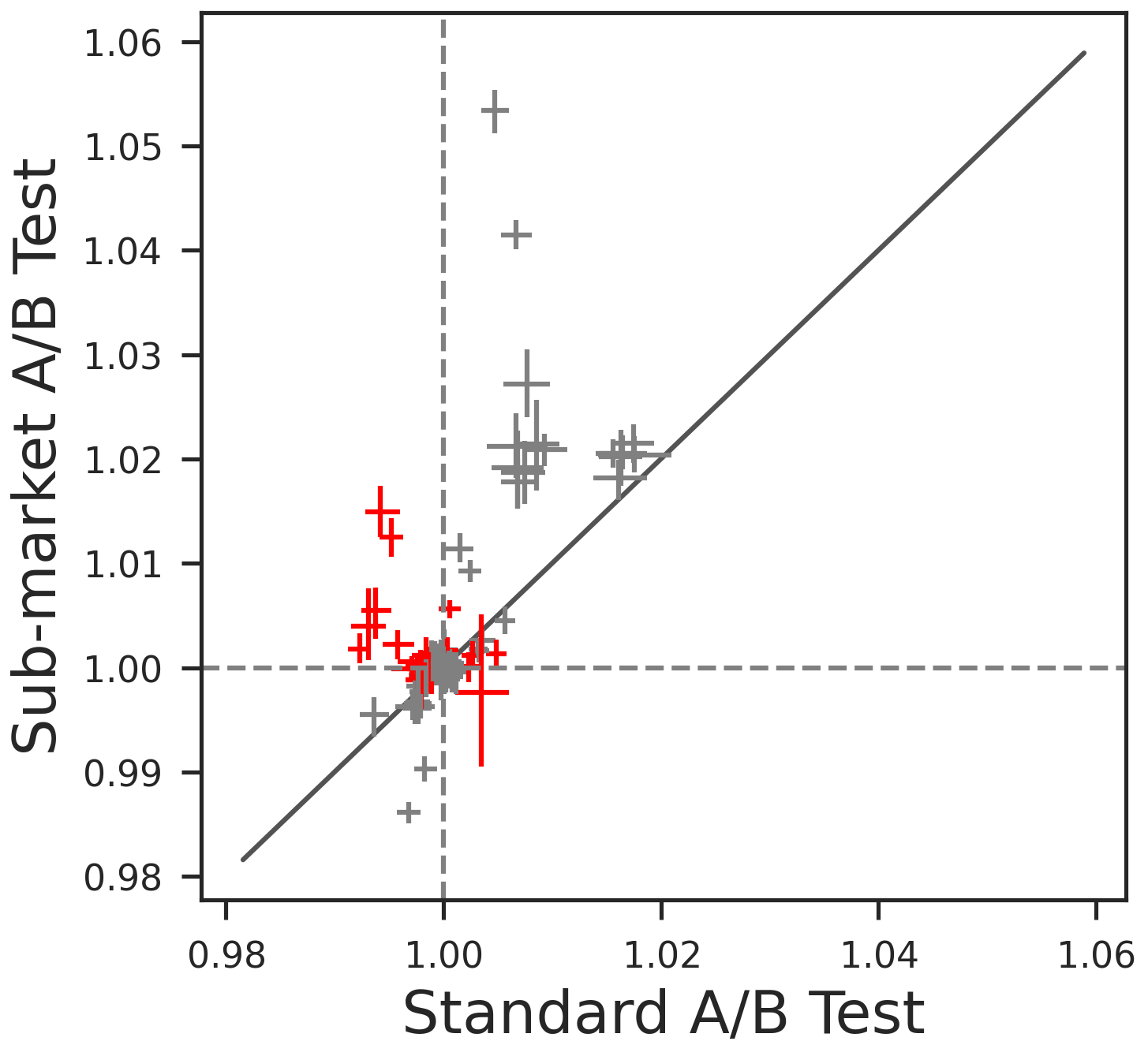}\qquad\qquad
    \includegraphics[width=0.34\linewidth]{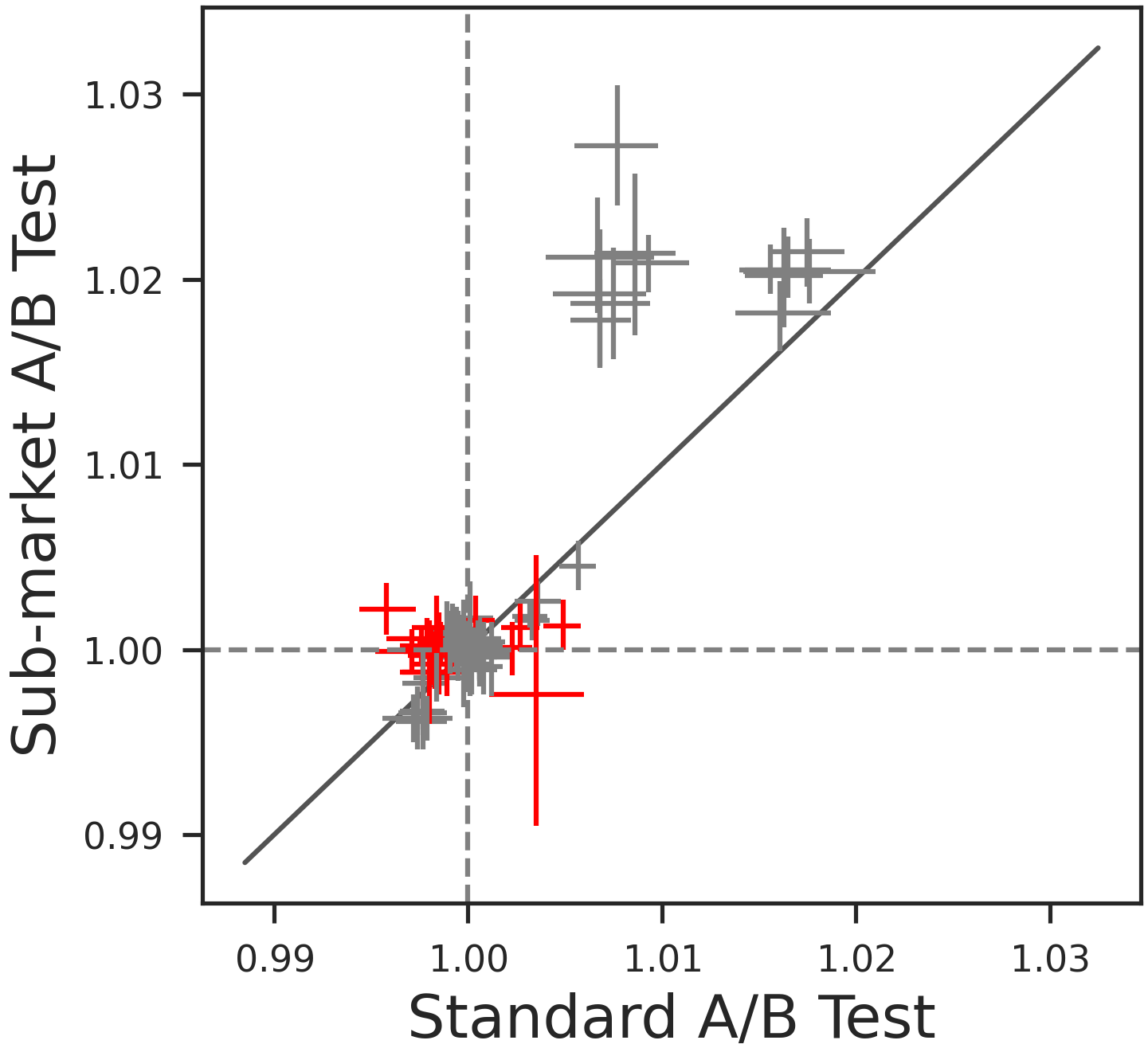}
	\label{fig:results}
	\caption{Parallel vs.\ standard budget-controlled A/B test, daily treatment effect. We denote neutral treatment effects with a value of $1.0$. Red crosses indicate instances of sign inconsistencies. Left are all datapoints, on the right, datapoints that fail a guardrail metric are removed.}
\end{figure}

Before the theoretical exposition, we consider a comparison of results for paired experiments between a parallel budget-controlled A/B test setup, and that of a traditional budget-split design; where the latter is considered the gold standard. \cref{fig:results} on the left shows comparisons of 99 experiments where the point estimate and CIs are plotted on the vertical axis for the parallel design, and on the horizontal axis for the budget-split design. The most important feature is whether the two experiments agree between (negative, neutral, positive), as a change would result in a launch reversal. The two experiment designs agree in $75\%$ of cases (at $90\%$ confidence level, hence the optimal agreement is $81.5\%$), which increases to $79\%$ after the introducing a guardrail metric, see \cref{fig:results} on the right. These results are quite satisfactory, but do point at the existence of remaining interference bias. In the remainder of this paper, correcting the interference bias is the main objective.

\paragraph{Contributions: }
\begin{enumerate}
    \item We formally define market interference in first-price auction markets using the first price pacing equilibrium (FPPE) framework \citep{conitzer2022pacing} (Sec. 3)
    \item We propose a debiased surrogate that eliminates the first-order bias of FPPE, and derive a plug-in estimator for the surrogate and establish its asymptotic normality. (Sec. 4)
    \item We run semi-synthetic experiments, confirming that the debiasing technique achieves the desired coverage properties. (Sec 5).
\end{enumerate}    

\subsection{Parallel A/B Testing in Practice}\label{ssc:inropractice}
In this section we describe the real-world problem of A/B testing with congestion that we wish to model, and our proposed solution of parallel A/B tests in carefully balanced submarkets. We start by describing the market environment. There is a set of $n$ advertisers, and each advertiser $i$ has a budget $b_i$. Whenever a user shows up on the platform an \emph{impression opportunity} occurs, and an auction is conducted in order to determine which ad will be shown to the user. Each advertiser $i$ has some stated value $v_{i}(\t)$ of being shown to a particular impression opportunity $\t$. The advertiser submits a bid which is determined based on $v_{i}(\t)$, as well as the expenditure of the advertiser so far. For example, in \emph{multiplicative pacing}~\cite{balseiro2017budget,conitzer2022multiplicative}, the platform adaptively learns a \emph{pacing multiplier} $\beta_i\in [0,1]$ such that the bid is formed as $ \beta_i v_{i}(\t)$. The budget-management system then adaptively controls $\beta_i$ over time, in order to ensure the correct rate of budget expenditure on behalf of the advertiser. We consider first-price auctions, which is the predominant way display advertising is sold online.

Now that we have discussed budget management, we describe the A/B testing problem. Suppose that a platform wants to run $K$ A/B tests, which may affect, e.g., the valuations that advertisers have for impression slots, revenue, etc. We construct the market-segmented experimental setup as follows: We first define a bipartite graph between advertisers and users based on targeting criteria. Subsequently, we cluster the advertisers into $K$ clusters, where the objective is to minimize the sum of weighted edges between clusters subject to traffic balancing constraints to make the resulting clusters as similar to the whole market as possible. The edge weight between a pair of advertisers is the number of impressions (or users) where they are both within the top-$k$ bids. The choice of $k$ is a parameter that must be chosen based on experience with the specific application setting.
If the clustering achieves a small objective function value, then each cluster is a mostly isolated submarket, in the sense that each user will mostly receive bids from advertisers in a single cluster.
Then, we run an A/B test within each of the $K$ submarkets.
Every user is randomly assigned to either ``A'' or ``B'' in each submarket. 
The main challenge is that while submarkets are relatively isolated, there is remaining interference from users who are targeted by advertisers from different submarkets, leading to a slightly different equilibrium. Our main contribution is to define a framework for analyzing such interference, and giving an estimator that removes the bias from these users. We survey related works in \cref{ssec:relatedwork}.

\section{Review of FPPE Theory}\label{sec:theoryoverview}

\textbf{Notation}.
For a measurable space $(\Theta, \diff \theta)$, we let $L^p$ (and $L^p_+$, resp.) denote the set of (nonnegative, resp.) $L^p$ functions on $\Theta$ w.r.t\ the base measure $\diff \theta $ for any $p\in [1, \infty]$ (including $p=\infty$).
Given $x \in L^\infty$ and $v \in L^1$, we let $\langle v, x \rangle = \int_\Theta v(\theta) x(\theta) \diff \theta$.
We treat all functions that agree on all but a measure-zero set as the same.
For a sequence of random variables $\{X_n\}$, we say $X_n = O_p(1)$ if for any $\epsilon > 0$ there exists a finite $M_\epsilon$ and a finite $N_\epsilon$ such that $\P(|X_n| > M_\epsilon) < \epsilon$ for all $n\geq N_\epsilon$. We say $X_n = o_p(1)$ if $X_n$ converges to zero in probability.
For a subset $\Theta'\subset \Theta$, let $1_{\Theta'} (\cd):\Theta \to \{0,1\}$ be the indicator function of $\Theta'$.
Convergence in distribution and probability is denoted by $\tod$ and $\toprob$.
Given a vector $a = [a_1, \dots, a_n]\tp$, let $\Diag(a)$ denote the diagonal matrix with $(i,i)$-th entry being $a_i$; sometimes we write $\Diag(a_i)$ when it is convenient to define each $a_i$ inline.
Let $\bm A\pinv$ denote the Moore--Penrose inverse of the matrix $\bm A$,
$\bm e_j$ the $j$-th unit vector,
and $[n]=\{1,\ldots,n\}$.

\begin{figure}[ht]
    \centering
    \begin{minipage}[c]{0.4\textwidth}
        \centering
        \includegraphics[width=\textwidth]{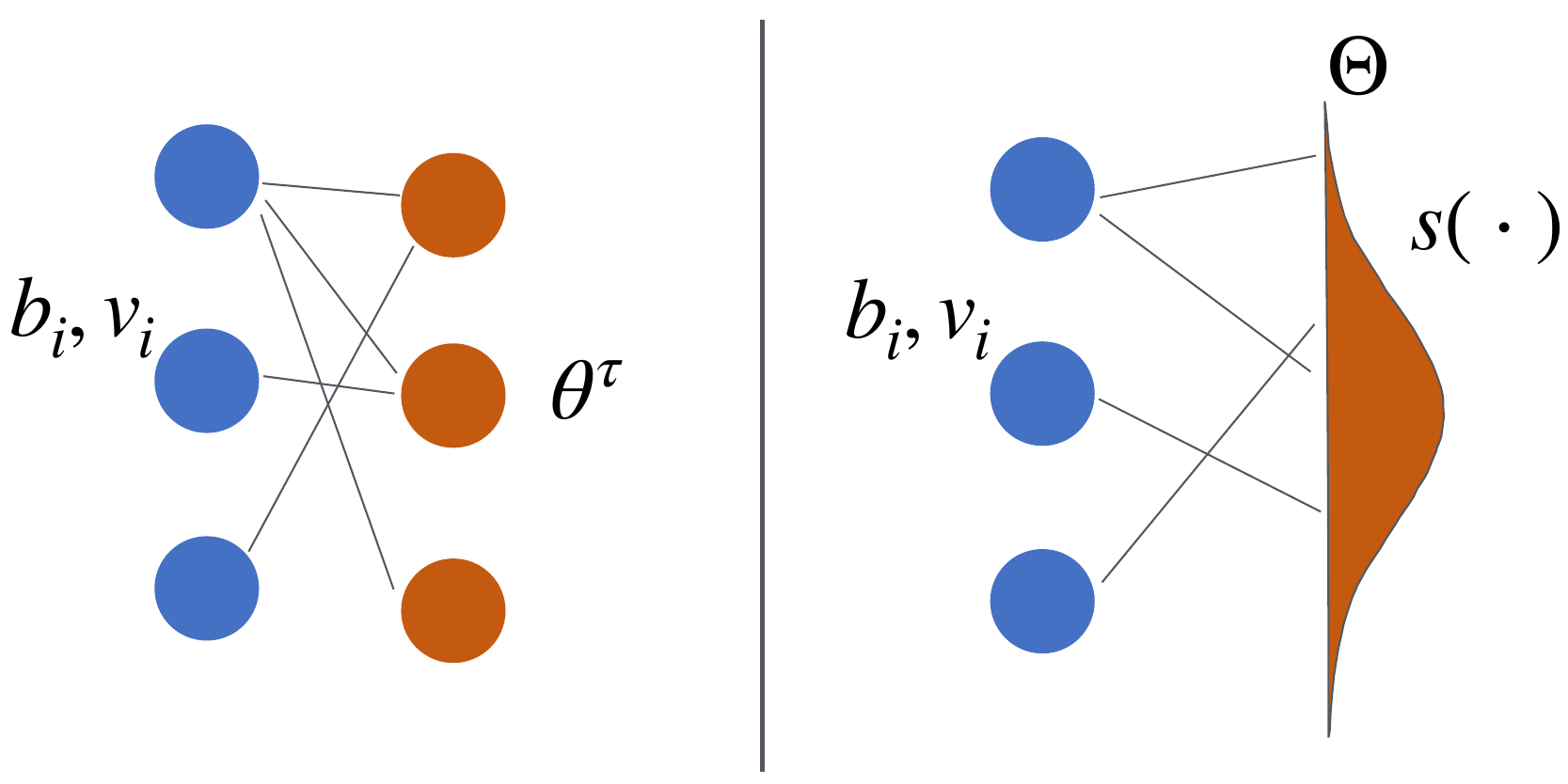}
    \end{minipage}
    \hfill
	\begin{minipage}[c]{0.5\textwidth}
        \centering
        \includegraphics[width=\textwidth]{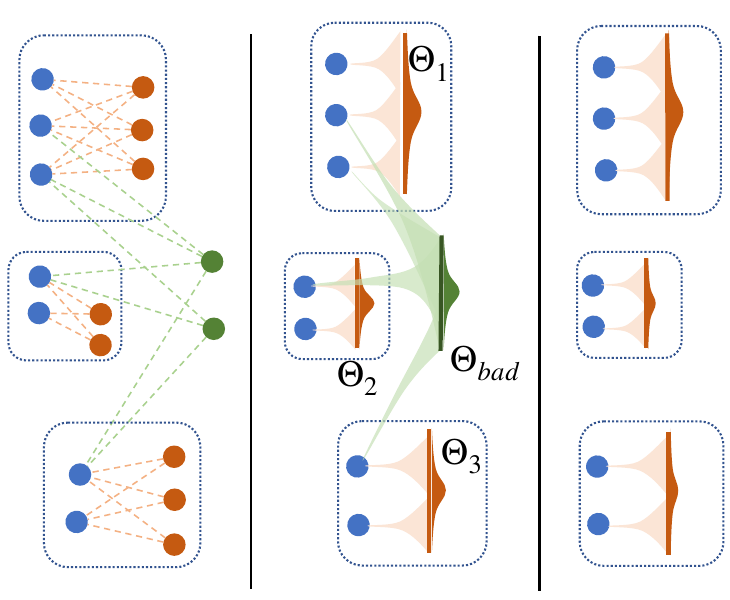}
    \end{minipage}
    \caption{
      \textbf{Left:} 
	  Finite FPPE (left) and limit FPPE (right). In a finite FPPE, there are a finite number of items; in a limit FPPE, the item set is a continuum.
      \textbf{Right:} The interference model --- 
      Left ($\widehat \cM\sa$): the observed market where interference is present among submarkets.
      Middle ($\cM\sa$): the limit market with interference from bad item set.
      Right ($\cM_0$): the limit market with perfectly separated submarkets.
      We use data from the left panel to make inferences about the market in the right panel.
    }
    \label{fig:combined}
\end{figure}

\noindent\textbf{Limit FPPE}.
We first introduce our notion of a limit market and two regularity conditions on the market, which models the underlying market structure that we sample from.
We have $n$ buyers and a possibly continuous set of items $\Theta$ with {an} integrating measure $\diff \theta$.
For example, one could take $\Theta = [0,1]$ and $\diff \theta  = $ the Lebesgue measure on $[0,1]$.
Each buyer has a \emph{budget} $b_i$; let $ \bm b = (b_1,\dots, b_n)$.
The \emph{valuation} for buyer $i$ is a function $v_i \in L^1_+$ such that buyer $i$ has valuation $v_i(\theta)$ for one unit of item $\theta\in \Theta$; let $\bm v: \Theta \to \Rn$, $\bm v(\theta) = [v_1(\theta),\dots, v_n(\theta)]$. We assume $\vbar = \max _i \sup_\theta  v_i(\theta)< \infty $.
The \emph{supplies} of items are given by a function  $ s \in L^\infty_+$, i.e., item $\theta\in \Theta$ has $s(\theta)$ units of supply.
	      Without loss of generality, we assume a unit total supply $\int_\Theta s \diff \theta = 1$.
	      Given $g:\Theta \to \R$, we let $\E[g] = \int g(\theta)s(\theta)\diff \theta$ and $\var[g] = \E[g\sq] - (\E[g])\sq$.
	      Given  $t$ i.i.d.\ draws $\{ \theta^1,\dots, \theta^t\}$
	      from $s$, let $P_t g(\cd) = \frac1t \sumtau g(\thetau)$.

Next we introduce the market equilibrium concept that is the foundation of our study. For that we leverage the \emph{first-price pacing equilibrium} (FPPE)~\citep{conitzer2022pacing}. FPPE models equilibrium outcomes under budget-management systems employed in several practical settings. Each buyer is assigned a \emph{pacing multiplier} $\beta_i$, which is used to control their budget expenditure. For each individual auction $\theta$, the buyer bids $\beta_iv_i(\theta)$, which can be seen as their adjusted valuation after factoring in their budget constraint (in practice, the valuation $v_i(\theta)$ may not be the buyer's true valuation, but instead their unscaled bid reported to the platform).
The goal of the budget management system is to achieve an equilibrium, meaning that it must ensure that buyers spend their budget exactly by appropriately choosing $\beta_i$. If the budget cannot be fully spent, then no pacing must occur (i.e., $\beta_i=1$).
Below we formally define the pacing equilibrium concept in the continuous setting (see \cref{fig:combined}, right), and give the finite version in the next section (\cref{fig:combined}, left).

\begin{defn}[Limit FPPE, \citet{gao2022infinite,conitzer2022pacing}]
	\label{def:limit_fppe}
	A limit FPPE, denoted as $\FPPE(\bm b , \bm v,s, \Theta)$, is the unique tuple $({\boldsymbol \beta}, p(\cd)) \in [0,1]^n \times L^1_+ (\Theta)$ such that there exist an allocation $x_i : \Theta \to [0,1]$ for all $i\in[n]$ such that,
	\begin{enumerate}[series = tobecont,itemjoin = \quad]
		\item (First-price)
		      Prices and allocations are determined by first price auctions:
		      for all items $\theta \in \Theta$, $p(\theta) = \max_i \betai v_i(\theta)$, and
		      only the highest bidders obtain items:
		      for all $i$ and $\theta$, $x_i(\theta) > 0$ implies $\betai \vithe =\max_k \beta_k v_k (\theta)$
		      \label{it:def:first_price}
		\item (Feasiblity, market clearing)
		      Budget are respected:
		      for all~$i$, $\int x_i(\theta) p(\theta) s(\theta)\diff \theta \leq b_i$.
		      There is no overselling:
		      for all $\theta$, $\sumiton x_i (\theta) \leq 1$.
		      \label{it:def:supply_and_budget}
		      Items with nonzero price are fully allocated:
		      for all $\theta$, $p(\theta) > 0$ implies $ \sumiton x_i(\theta) = 1$.
		\item (No unnecessary pacing) For all $i$, $\int x_i(\theta) p(\theta) s(\theta) \diff \theta < b_i$ implies $\betai = 1$.
		      \label{it:def:rev_max}
	\end{enumerate}
\end{defn}

Let $\betast$ and $\pst$ be the equilibrium pacing multipliers and prices.
Revenue in the limit FPPE is
$
	\REVst \defeq \int \pst(\theta) s(\theta)\diff \theta \; .
$
It measures the profitability of the auction platform.
The leftover budgets for buyers are denoted by $\deltasti \defeq b_i - \int \pst(\theta)s(\theta)\xsti(\theta)\diff \theta$.

The first two conditions simply describe the possible outcomes of a first-price auction system that uses pacing as the budget management strategy.
The last condition, no unnecessary pacing, ensures that we only scale down a buyer's bids in case their budget constraint is binding. FPPE has many nice properties, including that they are competitive equilibria and that they are revenue-maximizing among budget-feasible pacing multipliers~\citep{conitzer2022pacing}.

In a pacing auction market $\cM=\FPPE(\bm b , \bm v,s, \Theta)$ the following two regularity conditions are important for the study of its statistical properties.

\begin{defn}[\textsf{\scriptsize{SMO}}]
	\label{as:smoothness}
	We say the smoothness condition holds if
	the map ${\boldsymbol \beta} \mapsto \E_s[\max_i \betai \vithe ]$ is twice continuously differentiable in a neighborhood of $\betast$.
\end{defn}

\begin{defn}[\textsf{\scriptsize{SCS}}]
	\label{as:constraint_qualification}
	We say strict complementary slackness holds if,
	whenever a buyer is unpaced ($\betasti = 1$), then her leftover budget is strictly positive ($\deltasti  > 0$).
\end{defn}

The condition \nameref{as:smoothness} ensures that in the limit market, items that incur a tie are measure zero.
The condition \nameref{as:constraint_qualification} rules out degenerate buyers that spend their budget exactly at $\beta_i^*$. See \citet{liao2023statistical} and \citet{liao2023fisher} for an extensive discussion about these conditions.

\noindent\textbf{Finite FPPE}.
Next we introduce the finite FPPE, which models the auction data we observe in practice.
Let $\gamma = (\theta^1,\dots, \theta^t)$ be a sequence of items.
Assume each item has the same supply of $\sigma \in \Rp$ units.
A finite FPPE is a limit FPPE where the supply is a discrete measure supported on the observed items $\gamma$.
Let $\vitau = \vithetau$.

\begin{defn}[Finite FPPE, informal]    
    A finite FPPE, $\oFPPE(\bm b, \bm v,\sigma, \gamma)$, is a limit FPPE where the item set is the finite set of observed items $\gamma$. See \cref{sec:formal_def_finite_fppe} for the full definition.
\end{defn}


In \citet{liao2023statistical}, it is shown that if $\gamma$ consists of $t$ i.i.d.\ draws from distribution $s$, and one takes $\sigma = 1/t$, then the \pmr in  $\oFPPE(\bm b, \bm v,1/t, \gamma)$ converge to the \pmr in $\FPPE(\bm b , \bm v,s, \Theta)$ in probability.
Also, note that the FPPE $\oFPPE(t \bm b, \bm v, 1, \gamma)$ converges to the same limit FPPE $\FPPE(\bm b , \bm v,s, \Theta)$ because the pacing multipliers of a finite FPPE do not change when budgets and supplies are multiplied by the same scalar.

\noindent\textbf{The Eisenberg-Gale Program}.
Both the limit FPPE and the finite FPPE have convex program characterizations~\citep{chen2007note,conitzer2022pacing,gao2022infinite}.
We define the dual Eisenberg-Gale (EG) objective for a single item $\theta$ as
\begin{align}
	F (\t, \b) = f(\t,\b) -  \sumiton b_i \log \beta_i \;, \quad f(\t,\b)=\max_{i\in[n]} \beta_i v_i(\theta) \;.
\end{align}
The population and sample (dual) EG objectives are then defined as
\begin{align}
	\label{eq:def_pop_eg}
	H({\boldsymbol \beta}) = \E[F(\t,\b)] \;, \quad H_t({\boldsymbol \beta}) = \frac1t \sumtau F(\thetau, \b) \;.
\end{align}
We say $\cH$ is the Hessian of market $\cM$ when $\cH =\nabla_{\b\b} \sq \int F(\theta, {\boldsymbol \beta}) s \diff \theta | _{{\boldsymbol \beta} = \betast}$.

The equilibrium pacing multipliers $\betast$ in $\FPPE(\bm b , \bm v,s, \Theta)$ can be recovered through
the population dual EG program
\begin{align}\hspace{-.2cm}
	\label{eq:pop_deg}
	\betast = \argmin_{{\boldsymbol \beta} \in (0, 1]^{n}}
	H({\boldsymbol \beta})
	\;.
\end{align}
The pacing multiplier vector $\betast$ is the unique solution to \cref{eq:pop_deg}.
Let $\betagam$ be the equilibrium pacing multiplier in $\oFPPE(\bm b, \bm v,1/t, \gamma)$. Then $\betagam$ solves the sample analogue of \cref{eq:pop_deg}:
\begin{align}
	\label{eq:finite_deg}
	\betagam = \argmin_{\b\in(0,1]^n} H_t(\b)
	\;.
\end{align}

Let us briefly comment on the differential structure of $f$, since it plays a role in later sections.
The function $f(\b, \t)$ is a convex function of $\b$ and its subdifferential $\partial_\b f(\b,\t)$ is the convex hull of $\{ v_i \bm e_i: \text{index $i$ such that $\betai\vithe = \max_k \betak v_k(\t)$}  \}$, with $\bm e_i$ being the base vector in $\Rn$.
When $\max_i \betai v_i(\t)$ is attained by a unique $i^*$, the function $f(\cdot,\theta)$ is differentiable. In that case, all entries of $\nabla_\b f(\b,\t)$ are zero expect that the $i^*$-th entry is filled with the value $v_{i^*}(\t)$.

\section{Interference as Contamination}\label{sec:setup}

In this section, we discuss how to estimate market equilibria when there is \emph{contamination} in the supply, meaning that items are generated from a mixture of two distributions, when in reality we wish to estimate equilibrium quantities from one of the two distributions. Then, we show that interference from other markets can be viewed as a form of contamination, and so the problem of removing interference bias can be analyzed via our contamination framework.

\subsection{FPPE with Contaminated Supply}

We assume that we are in the same FPPE setting as before:
there are $n$ buyers, each with budget $b_i$,
and an item set $\Theta$ which is now partitioned in to $\Theta_\bad$ and $\Theta_\good$. However, now we assume that the supply $s$ is contaminated by $s'$, another supply distribution.
We define the $\alpha$-contaminated market as $\cM_\alpha = \FPPE(\bm b, \bm v, s \sa , \Theta)$
and the uncontaminated market as $\cM_0 =  \FPPE( \bm b, \bm v,  s, \Theta)$,
where $s\sa = \alpha s' + (1-\alpha) s$, distribution $s'$ is supported on $\Theta_\bad$, and $s$ on $\Theta_\good$.
Our goal is to perform inference about FPPE properties in the limit FPPE with the supply $s$. However, we are given access to finite FPPEs sampled from $s_\alpha$ instead.
In particular, let $\gamma$ be $t$ i.i.d.\ draws from $s\sa$ and let
	$\widehat \cM\sa = \oFPPE(\bm b, \bm v,1/t, \gamma)$.
We assume $\alpha$ is known throughout the paper.
In practice, this can often be estimated from historical data; in the parallel A/B test setting, this can be estimated directly from the sampled set of items, since we know whether an item is drawn from $s$ or $s'$.
Let $\betast$ and $\betasta$ be the limit pacing multipliers in $\cM_0$ and $\cM\sa$, respectively. Let $\betagama$ be the pacing multipliers in the sampled market $\widehat \cM\sa$ and let $H_{\alpha,t}(\b) = \frac1t \sumtau F(\thetau, \b)$ be the sample EG objective.

If we wanted to make inferences about $\cM\sa$ then we could use existing statistical inference theory on how to use data in a finite FPPE $\widehat \cM\sa$ to make inferences about the limit FPPE $\cM\sa$~\citep{liao2023statistical}. However, the supply contamination prevents the application of these tools to our problem.

Our central research question is then on how to use data in the finite contaminated market $\widehat \cM\sa$ to make inferences about the uncontaminated limit market $\cM_0$.

In \cref{sec:results} we propose an estimator for this problem and derive its properties. The results there apply to general item space $\Theta$ and supplies $s$ and $s'$. By imposing structure on $\Theta, s$ and $s'$, we show that the contamination model captures the interference among FPPEs. 

\subsection{Application: Modeling Interference among FPPEs} \label{sec:inf_as_con}
Now we show how the contamination model from the previous section can be used to model interference.
Consider $K$ separate auction markets, which together form a global market.
In the global market there are $n$ buyers, each with budget $b_i$,
and an item set $\Theta$, partitioned into $\Theta_\good$ and $\Theta_\bad$.
Let $C_1, \dots, C_K$ be a partition of the buyers,
$\Theta_1, \dots, \Theta_K$ be a partition of the good item set $\Theta_\good$,
and $s_1,\dots s_K$ be a set of supply functions, supported on $\Theta_1,\dots, \Theta_K$ respectively.
The $k$-th submarket
consists of buyers in $C_k $, the item set $\Theta_k$ and supply $s_k$.
Let $s = \frac1K \sum_k s_k$ be the average mixture and $s'$ be a supply supported on $\Theta_\bad$.
Let the contaminated supply be $s\sa = \alpha s' + (1-\alpha) s$.

By imposing structure on
$\Theta_\good$ and $\Theta_\bad$
the contamination model can capture interference among auction markets.
We assume that submarkets are separated, which models the ideal case where there is no interference.
A buyer $i \in C_k$ is only interested in items from the submarket he belongs to: $v_i(\theta) = 0$ for $\theta \in \Theta_{k'}, k'\neq k$.
Next, we let $\Theta_\bad$ represent items
that cause outbound edges from submarkets; see the green edges in \cref{fig:combined} left panel. An item is referred to as {\em bad} if it has positive values for buyers from at least two different submarkets. Formally,
$\theta \in \Theta_\bad$ if there exist $i \in C_k$, $j\in C_{k'}$, $k\neq k'$, such that $v_i(\theta) > 0$ and $v_j (\theta) > 0$.
Combining these assumptions, we have that a buyer $i$ from submarket $k$ has positive values \emph{only} for items from the sets $\Theta_k$ and possibly $\Theta_\bad$.
Now we have fully specified a contaminated market setup: we wish to make inferences on the market consisting of only $\Theta_\good$ (which is really $K$ fully separate submarkets), but we observe an actual market containing items from $\Theta_\good \cup \Theta_\bad$.
With this setup, we can use the results developed in the following section to model interference in parallel submarkets.

In \cref{fig:combined} we present an example of interference among $K=3$ submarkets.
The market of interest is the perfectly separated market (right).
This is because, in parallel A/B testing, submarkets are explicitly created such that each submarket resembles the global market. Then when a submarket receives a treatment, the observed quantities in that submarket, such as revenues and social welfare, are considered surrogates for the treatment effect in the global market.
However, in practice we only observe the interfered finite market (left), which converges to the interfered limit market (middle). In \cref{sec:ABTest} we show how to analyze parallel A/B testing using this framework.

\section{A Debiased Estimator and Its Properties}\label{sec:results}
This section develops a methodology for making inferences about the uncontaminated limit FPPE.
Since the interference setting is a special case of the contamination setting, we develop theories for the latter.
We introduce a surrogate for pacing multipliers, based on the notion of directional derivatives, and establish its debiasing property in~\cref{sec:def_surr}. Then, we focus on estimating this surrogate quantity in~\cref{sec:estimator}, and develop asymptotic normality results in~\cref{sec:normality}.
Secondly, we consider estimating revenue, which can be thought of as a smooth function of pacing multipliers. We discuss debiased revenue estimation and inference based on our pacing multiplier results in \cref{sec:res_rev}.



\subsection{A Debiased Surrogate for Pacing Multipliers} \label{sec:def_surr}

If we view $\betasta$ as a function of the level of contamination $\alpha$, then one can imagine that under sufficient regularity conditions, the pacing multipliers in the perfectly separated market, $ \betast_0$, can be approximated by some form of Taylor expansion of $\alpha \mapsto \betasta$ at $\alpha$.
This can be made rigorous by the notion of directional derivatives. We define
\begin{align}
	\diff \betast(\alpha )= \lim_{ \epsilon \to 0^+} \frac{  \betast_{\alpha - \epsilon} - \betasta }{\epsilon}
	\;.
\end{align}
if the limit exists.
We will show in \cref{thm:bias} that $\betasta + \alpha \diff \betast(\alpha)$ serves as a good approximation to $\betast = \betast_0$.

Thanks to the convex program characterization of FPPE, the directional derivative $\biasformula$ has a closed-form expression under certain regularity conditions (the conditions are given in \cref{thm:bias}; the full proof is given in the appendix). We need a few notations for this expression.
Define
$
	{\boldsymbol \delta} \sa
	=  \int \nabla f (\theta, \betasta) (s - s') \diff \theta .
$
Let $\cH\sa = \nabla_{\b\b}\sq \int F(\t,\betasta) s\sa\diff\t$ be the Hessian matrix in the market $\cM\sa $ and $  {\boldsymbol P} \sa = \Diag( 1 (\betast_{\alpha, i} < 1) )$. Then, under
the regularity conditions given in \cref{thm:bias} below,
\begin{align} \label{eq:defbias}
	\biasformula  = - ({\boldsymbol P} \sa \cH\sa {\boldsymbol P} \sa )\pinv {\boldsymbol \delta} \sa
	\;.
\end{align}
We present a heuristic derivation in \cref{sec:heuristic}.
Given the closed-form expression of $\biasformula$, we define the following debiased pacing multiplier
\begin{align}
	\label{eq:def_debiased_surr}
	\betastt  = \betasta + \alpha \cdot (- ( \boldsymbol{P} \sa \cH\sa {\boldsymbol P} \sa )\pinv {\boldsymbol \delta} \sa)
	\;.
\end{align}

\begin{theorem}[Analysis of Bias]
	\label{thm:bias}
	Suppose that in the market $\cM_0$ conditions \nameref{as:smoothness} and \nameref{as:constraint_qualification} hold, and assume that
	$\b \mapsto \nabla \sq \int F(\t,\b)s'\diff\t$ is twice continuously differentiable in a neighborhood of $\betast$.
	Then the directional derivative $\diff \betast(\cd)$ is well-defined in a neighborhood of zero, and is given by \cref{eq:defbias}.
	Moreover, as $\alpha \downarrow 0$,
	\begin{align*}
		\| \betastt - \betast \|_2 = o (\alpha ) \;.
	\end{align*}
\end{theorem}
The proof is given in \cref{sec:proof:thm:bias}.
\cref{thm:bias} indicates that the debiased surrogate $\betastt$ removes first-order bias
caused by contamination.
The limit pacing multipliers $\betasta$ of the contaminated market $\cM\sa$ will have bias of order $\betasta - \betast = \Theta(\alpha)$. In contrast, \cref{thm:bias} shows that
the debiased surrogate only incurs a bias of order $o(\alpha)$.

\subsection{The Estimator} \label{sec:estimator}
In this section we introduce a plug-in estimator for the debiased surrogate $\betastt$ and introduce a consistency theorem.
The next section discusses constructing confidence intervals.

To estimate $\biasformula$ in \cref{eq:defbias} we need estimates of its three components: the Hessian $\cH\sa = \nabla_{\b\b}\sq \int F(\t,\betasta) s\sa\diff\t$, the diagonal matrix ${\boldsymbol P} \sa$ and the vector ${\boldsymbol \delta} \sa = \int \nabla f (\t, \betasta) (s-s')\diff\t$.

\noindent\textbf{The Hessian}. For simplicity in our theoretical results, we will simply assume a generic Hessian estimator $\widehat \cH \sa$ such that for some $\eta_t \downarrow 0$ we have
	      $\widehat \cH\sa - \cH\sa = O_p( \eta_t) $.
	      \citet{hong2015extremum} discuss the estimation of the derivative in detail.
	      Different kinds of statistical guarantees require different rate conditions on $\eta_t$; see \cref{thm:betaconsistency,thm:betaclt}.
	      We then introduce two Hessian estimators: one is applicable for general FPPE, while the other requires an extra market regularity condition.
	      The first Hessian estimator is the finite difference method.
	      Let $\bm e_i, \bm e_j$ be basis vectors and $\varepsilon_t$ be a step-size.
	      Then the estimator is
	      $$
		      \widehat \cH\sa[i,j] = [H_{\alpha, t}(\betagam_{++})-H_{\alpha, t}(\betagam_{+-})
		      -H_{\alpha, t}(\betagam_{-+})
		      +H_{\alpha, t}(\betagam_{--})] /
		      (4 \varepsilon_t^2)  \;,
	      $$
	      where $\betagam_{\pm \pm} \defeq  \betagama \pm \bm e_i \varepsilon_t \pm \bm e_j \varepsilon_t$, and $H_{\alpha, t}({\boldsymbol \beta}) = \frac1t \sumtau F(\thetau, {\boldsymbol \beta})$, with $\{\thetau\}_\tau$ being the items in $\widehat \cM \sa$.
	      In practice, a diagonal approximation of the Hessian suffices.
	      The second method relies on an additional regularity condition, in which case we derive a simplified formula for the Hessian, thereby enabling a simpler estimation procedure (see \cref{thm:betaclt}).

\noindent \textbf{The vector ${\boldsymbol \delta} \sa$.}
	      Let  $g$ be the Radon-Nikodym ratio $g (\theta) = (\diff  (s - s') / \diff s \sa ) (\theta) =  \frac{1}{1-\alpha } 1_{\Theta _\good}(\theta) -\frac{1}{\alpha } 1_{\Theta_\bad} (\theta)$.
	      With the ratio $g$, the true vector ${\boldsymbol \delta} \sa$ can be written as ${\boldsymbol \delta} \sa = \int g(\theta) \nabla f(\theta, \betasta)  s\sa(\theta)\diff\theta$, which is easy to estimate given i.i.d.\ draws from $s\sa$.
	      In particular, our estimator is then
	      $
		      \widehat {\boldsymbol \delta} \sa = \frac1t \sumtau g(\thetau) \mutau \; .
	      $
	      Here $\mutau = [x_1^\tau v_1^\tau, \dots, x_n^\tau v_n^\tau]\tp $
	      is a subgradient of $f(\thetau, \betagama)$ w.r.t.\ ${\boldsymbol \beta}$.

\noindent \textbf{The diagonal matrix ${\boldsymbol P} \sa$.} Recall ${\boldsymbol P} \sa = \Diag( 1(\betastai < 1))$. So a natural estimator is
	      $        \widehat {\boldsymbol P} \sa = \Diag ( 1(\betagamai < 1 - \iota_t) )$,
	      where the slackness $\iota_t \asymp \frac{1}{\sqrt t}$.

With all three components estimated, using  $\betagama$ is the \pmr in the market $\widehat \cM\sa$ we define the plug-in estimator for $\biasformula$ in \cref{eq:defbias} as
\begin{align}
	\label{eq:def_estimator}
	\widehat{\boldsymbol \beta} = \betagama -  \alpha  \,\cd \,  (\widehat {\boldsymbol P} \sa \widehat \cH\sa \widehat {\boldsymbol P} \sa)\pinv \widehat {\boldsymbol \delta} \sa
	\;.
\end{align}

\begin{theorem}[Consistency] \label{thm:betaconsistency}
	Suppose that in the market $\cM\sa$ conditions \nameref{as:smoothness} and \nameref{as:constraint_qualification} hold. 
	If the Hessian estimation error satisfies $\eta_t = o(1)$, then $\widehat {\boldsymbol \beta} \toprob \betastt$. The proof is in \cref{sec:proof:thm:betaconsistency}.
\end{theorem}

\subsection{Asymptotic Normality and Inference} \label{sec:normality}
We present two asymptotic normality results.
In the first result, we require a
stronger condition on the Hessian error rate $\eta_t$.
In particular,  as will be shown in \cref{thm:betaclt}, the rate condition $\eta_t = o(1/\sqrt t)$ is sufficient for normality.
One could use a separate large historical dataset to obtain a good estimate of the Hessian matrix.
In the second result, we
impose an additional condition on market structure which simplifies the Hessian expression and facilitates efficient Hessian estimation.
To describe the additional market structure, we define the gap between the highest and the second-highest bid for an item $\theta$ under pacing ${\boldsymbol \beta}$ by
$
	\bidgap({\boldsymbol \beta},\theta) \defeq \max \{\betai \vithe  \} - \operatorname{secondmax}\{ \betai \vithe \} \;,
$
where $\operatorname{secondmax}$ is the second-highest entry
potentially equal to the highest; e.g., $\operatornamewithlimits{secondmax}([1,1,2]) = 1$.
When there is a tie for an item $\theta$ under pacing ${\boldsymbol \beta}$, we have $\bidgap({\boldsymbol \beta},\theta) = 0$.
When there is no tie for an item $\theta$, the gap $\bidgap({\boldsymbol \beta},\theta)$ is strictly positive.

For any $g: \Theta \to \R$, let $\E\sa [g] = \int g s\sa\diff\theta$ and $\cov\sa (g) = \E\sa[ (g - \E\sa[g]) (g - \E\sa[g]) \tp]$.
Recall $\betastt$ is the debiased surrogate in \cref{eq:def_debiased_surr} and $\widehat {\boldsymbol \beta}$ is its estimator defined in \cref{eq:def_estimator}.

We need to introduce a few more notations to describe the normality results. First, let
$
	\bm d_\alpha = - ( {\boldsymbol P} \sa \cH\sa {\boldsymbol P} \sa )\pinv \nabla f(\cd, \betasta) \;$.
As mentioned previously, the pacing multipliers in the contaminated market converge to the limit counterpart  and have the representation
\begin{align*}
	\sqrt t (\betagama - \betasta) = \frac{1}{\sqrt t} \sumtau ( \bm d_\alpha(\thetau) - \E\sa[\bm d_\alpha(\thetau) ] )+ o_p(1)
	\;.
\end{align*}
In the statistics literature, the function $\bm d_\alpha(\cd)-\E[\bm d_\alpha]$ is called the influence function \citep{van2000asymptotic}.
For our debiased estimators, we need the following (uncentered) influence functions.
\begin{align*}
	  \IF_1(\t) = \big( \tfrac{1 }{1-\alpha} 1_{\Theta_\good} (\t) \big)\bm d_\alpha (\t) \;,
	 \; \IF_2 (\t) = \IF_1 (\theta) - 2\alpha \Diag( \betastai  \delta_{\alpha ,i} b_i \inv )  \bm d_\alpha (\theta) \;.
\end{align*}

\begin{theorem}\label{thm:betaclt}
	Let
	\nameref{as:constraint_qualification} and \nameref{as:smoothness} hold in $\cM\sa$.	
    
    1. \textbf{Asymptotic Normality in a General Market.}    It holds that
	$
		\widehat {\boldsymbol \beta} - \betastt  = {\bm z_t} + O_p(\eta_t) + o_p(\tfrac{1}{\sqrt t }),
	$
	where
	$
		\sqrt t \bm z_t \tod \cN(\bm 0, {\bm \Sigma}_1)
    $
	with ${\bm \Sigma}_1 = \cov\sa(\IF_1)$, and $\eta_t$ is the Hessian estimation error.
    
	2. \textbf{Asymptotic Normality under a Bid Gap Condition.} If in addition $\E\sa [1 / \bidgap(\betasta, \theta)] < \infty$, then $\cH\sa = \Diag( b_i / (\betastai)\sq )$. Suppose we estimate $\cH\sa$ with $ \widehat \cH \sa = \Diag(b_i / (\betagamai)\sq )$. Then
	$
		\sqrt t (\widehat {\boldsymbol \beta} - \betastt) \tod \cN(0, {\bm \Sigma}_2)
	$
	where ${\bm \Sigma} _ 2 = \cov\sa (\IF_2)$.
	The proof is in \cref{sec:proof:thm:betaclt}.
\end{theorem}
\cref{thm:betaclt} part 1 shows how the error of the Hessian estimate affects the distribution of the estimator $\widehat {\boldsymbol \beta}$.
If $\eta_t = o(1/\sqrt t)$, then
the decomposition becomes $\sqrt {t} (\widehat{\boldsymbol \beta} - \betastt) = \sqrt {t} \bm z_t + o_p(1)$, implying asymptotic normality, i.e.
$\sqrt t (\widehat {\boldsymbol \beta} - \betastt) \tod  \cN(\bm 0, {\bm \Sigma}_1)$, in which case one can construct an ellipsoidal confidence region for $\betastt$.
\cref{thm:betaclt} part 2 shows direct
asymptotic normality under the extra condition, with a simpler Hessian estimator that avoids finite differences.

To perform inference, we need to construct a consistent estimate of the covariance matrix.
Now we describe a plug-in estimate of ${\bm \Sigma}_1$. Let the estimator $\widehat{\IF_1}^\tau$ be
$
\widehat{\IF_1}^\tau = -\big( \tfrac{ 1 }{1-\alpha} 1_{\Theta_\good} (\thetau) \big)(\widehat {\boldsymbol P} \sa \widehat \cH\sa \widehat {\boldsymbol P} \sa) \pinv \mutau
	\;,
$
where $\mutau = [x_1^\tau v_1^\tau, \dots, x_n^\tau v_n^\tau]\tp $, and $\widehat{{\boldsymbol P} \sa}, \widehat \cH\sa$ have been defined in \cref{sec:estimator}.
The plug-in estimator is
$
	\widehat {\bm \Sigma}_1 = \frac1t \sumtau (  \widehat{\IF_1}^\tau - \overline{\IF_1}) ( \widehat{\IF_1}^\tau - \overline{\IF_1}) \tp
$
with $\overline{\IF_1} = \frac1t \sumtau \widehat{\IF_1}^\tau$.
By similar arguments as in \citet{liao2023fisher}, the plug-in estimates of ${\bm \Sigma}_1$ and ${\bm \Sigma}_2$ are consistent. \cref{alg}  summarizes the debiasing procedure.

In \cref{sec:res_rev} we present a similar debiased estimator for revenue and its bias and variance properties. In \cref{sec:ABTest} we specialize the debiased estimator to parallel budget-controlled A/B testing.

\section{Semi-synthetic experiment} \label{sec:exp}
To evaluate our proposed framework and debiased estimator, we run semi-synthetic simulations to check if the proposed estimator for beta and revenue are indeed less biased, and we test the coverage of the proposed estimator. Fully synthetic experiments are presented in~\cref{ssec: simulation}. 

In the semi-synthetic experiments, we simulate 40 buyers and 10000 \emph{good} items in two submarkets, with a varying number of \emph{bad} items (up to 5000) in order to study the effect of the contamination parameter $\alpha$. For each $\alpha$, we randomly sample a budget for each buyer, and compute ${\boldsymbol \beta}^*$ and $\REV^*$ from the limit pure market $\cM_0$ with a value function in each submarket. Both the budget and values are sampled from historical bidding data, making the budget and value distributions heavy-tailed as in the real-world applications. More specifically, we first sample a certain number of auctions. For each auction, we sample a given number of advertisers with their per-impression bids. Advertisers that are sampled across different auctions are treated as the same buyers and their budgets are determined by aggregating their values over auctions up to a scalar to calibrated
to get the percentage of budget-constrained buyers equal to what was observed in the real-world auction market, along the same lines as the experiments of \citet{conitzer2022multiplicative}.

\begin{figure*}
	\centering
	\includegraphics[width=0.8\textwidth]{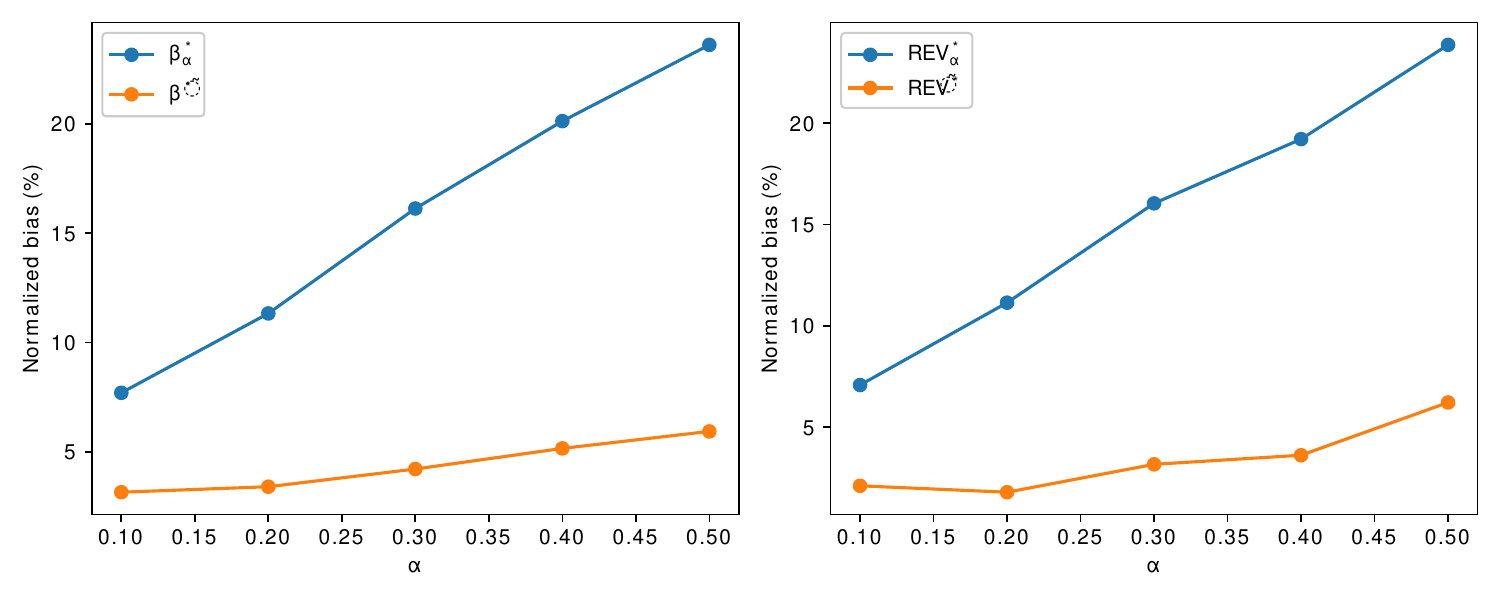}
	\caption{Normalized bias (in percent of true value) as a function of $\alpha$ in semi-synthetic experiments. $\widetilde{\boldsymbol \beta}^*$ and $\widetilde\REV^*$ are the debiased surrogates for pacing multiplier and revenue in the limit market with interference $\cM\sa$. }
	\label{fig:semi_sim_bias}
\end{figure*}

To check if the debiased surrogate reduces bias, we compute ${\boldsymbol \beta}^*\sa$ and $\REV^*\sa$ from the limit market with interference $\cM\sa$ and their surrogates, $\widetilde{\boldsymbol \beta}^*$ and $\widetilde\REV^*$ ($\widetilde{\boldsymbol \beta}^*$ is defined \cref{eq:def_debiased_surr}, $\widetilde\REV^*$ is defined in \cref{sec:res_rev}. We look at the normalized bias for the surrogate, defined as $\| \betastt - \betast \|_2 / \| \betast \|_2$ for pacing multipliers and as $|\widetilde\REV^* / \REV^* - 1|$ for revenue, and similarly defined for the limit quantities. \cref{fig:semi_sim_bias} shows the normalized bias curves as a function of $\alpha$.
The magnitude of the bias increases with $\alpha$, for both the variables in the limit market with interference $\cM\sa$ and their debiased surrogates. The bias of the debiased surrogates is indeed much smaller than the contaminated limit quantities.

Next, we check the coverage of the proposed variance estimator. For each $\alpha$ and each budget sample, we run 100 simulations in the following way: We sample items (or their values for each buyer) considering two submarkets and bad items. We then run the finite FPPE with bad items and obtain a baseline estimate for pacing multiplier and revenue without applying the debiasing procedure. Then, we apply the debiasing procedure to compute the debiased estimates. For each simulation, we check if the debiased surrogate is within the confidence interval of the debiased estimator. Finally, we aggregate them to compute the estimated coverage of the estimator. The results for both pacing multiplier and revenue are shown in \cref{tab:semi_synthetic_revenue_coverage}. For the coverage of $\widehat{\boldsymbol \beta}$, we first compute the coverage of each component and report only the average in the table. For revenue, we construct the CI using the two approaches as mentioned in~\cref{sec:res_rev}: one based on \cref{eq: revenue_ci} and the other using parametric bootstrap based on the estimated asymptotic distribution of $\widehat{\boldsymbol \beta}$ (with "(b)" in the column names). 

\begin{figure}
	\centering
	\includegraphics[width=0.33\textwidth]{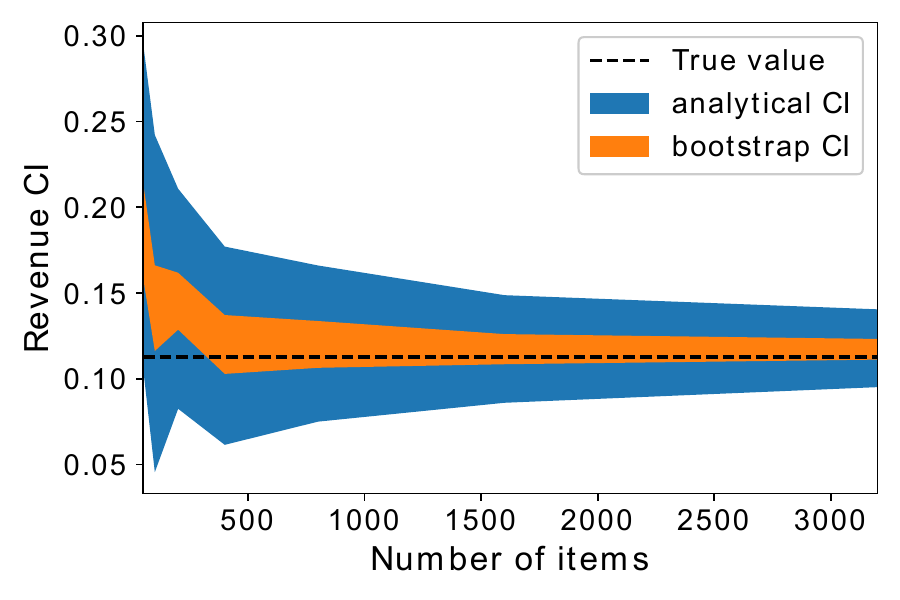}
	\caption{Revenue confidence intervals as a function of the number of items in semi-synthetic experiments. The analytic CI comes from \cref{eq: revenue_ci}. The true value is the debiased surrogates for revenue in the limit market with interference $\cM\sa$. 
 }
	\label{fig:semi_sim_convergence}
\end{figure}

Firstly, \cref{fig:semi_sim_convergence} shows that both CIs converges to the true value from the limit market with interfrence $\cM\sa$ as the number of items goes to infinity. Then, in ~\cref{tab:semi_synthetic_revenue_coverage}, we show that the coverage for $\widehat{\boldsymbol \beta}$ is slightly smaller than the nominal level (95\%), as well as the coverage of the bootstrap CI of revenue. The under-coverage for $\widehat{\boldsymbol \beta}$ is mainly driven by the under-estimation of the variance of $\widehat{\boldsymbol \beta}$, while the under-coverage of the bootstrap CI for revenue can also be partially attributed to the higher dimensionality (with 40 buyers), making the bootstrap resampling harder to explore the whole space. 

Although the proposed variance estimator has good asymptotic properties, the results from our synthetic experiments suggest that it can perform badly, in either direction, for finite markets. Constructing more accurate variance estimators for our debiased estimator in finite settings would certainly mitigate the over- or under-coverage issues that we observe here and deserve more future research. One promising alternative is to construct the CI for $\widehat{\boldsymbol \beta}$ and $\widehat\REV$ by directly bootstrapping the observed value matrix, though this might work best for independent valuations across buyers.

\begin{table}[!ht]
	\centering\small
	\begin{tabular}{c c c c c c c}
		\toprule
		           & \multicolumn{1}{c}{$\widehat{\boldsymbol \beta}$} & \multicolumn{4}{c}{$\widehat\REV$}                                            \\
		\cmidrule(lr){2-2} \cmidrule(lr){3-6}
		$\alpha$   & coverage                            & CI width                           & CI width (b) & coverage & coverage (b) & \\
		\midrule
		1000/11000 & 0.877                               & 0.244                              & 0.044        & 1.0      & 0.95           \\
		2000/12000 & 0.849                               & 0.225                              & 0.043        & 1.0      & 0.87           \\
		3000/13000 & 0.852                               & 0.210                              & 0.041        & 1.0      & 0.81           \\
		4000/14000 & 0.828                               & 0.200                              & 0.040        & 1.0      & 0.90           \\
		5000/15000 & 0.826                               & 0.191                              & 0.039        & 1.0      & 0.90           \\
		\bottomrule
	\end{tabular}
	\caption{Coverage of $\widehat{\boldsymbol \beta}$ and revenue estimates in the semi-synthetic experiments. All quantities are averaged over 100 simulations for each $\alpha$ (the ratio of the number of bad items and the total items). The coverage of $\widehat{\boldsymbol \beta}$ is averaged over all components of $\widehat{\boldsymbol \beta}$. For revenue estimates, the columns with "(b)" represent the quantities from the bootstrap CI, while the columns without "(b)" are for the CI from \cref{eq: revenue_ci}. The CI widths are normalized by the revenue from the limit market $\cM\sa$.}
	\label{tab:semi_synthetic_revenue_coverage}
\end{table}

\section{Conclusion} \label{sec:conclusion}
We have proposed a practical experimental design for performing concurrent A/B tests in large-scale ad auction markets, using a submarket clustering approach, and showed that in production experiments, this submarket clustering approach leads to strong sign consistency performance, as compared to A/B testing on the full market, while allowing significantly-higher A/B test throughput.
In order to model the potential for interference between submarket A/B tests, we introduced a theoretical model of statistical inference in first-price pacing equilibrium problems, under settings with supply contamination. We showed how one can perform statistical inference in such a setting using a debiased estimation procedure, and studied the statistical properties of this procedure.
We then showed how our model of statistical inference in FPPE with contamination can be used to model the submarket clustering parallel A/B test design, and gave theoretical performance guarantees.
Finally, we presented numerical experiments on fully synthetic and semi-synthetic data derived from Meta ad auctions. The experiments showed that our proposed debiased estimator achieves smaller biases and its statistical coverage on realistic data is generally in line with the predictions from our theory.

\newpage 
\section*{Acknowledgements}

We would like to thank the anonymous reviewers for their useful comments.
Christian Kroer and Luofeng Liao were supported by the Office of Naval Research awards N00014-22-1-2530 and N00014-23-1-2374, and the National Science Foundation awards IIS-2147361 and IIS-2238960.


%% file: app_theory.tex
\section{Technical Lemmas}

\subsection{A perturbation result for constrained stochastic programs}
We introduce a theorem from \citet{shapiro1990differential}.
Theorem~1 in that paper handles the case where the constraints are also defined as expectations of some random functions, and constraints are also perturbed in the analysis. For simplicity we specialize the theorem to deterministic constraints and do not perturb constraints. Note that Theorem~7.27 from \citet{shapiro2021lectures} can also be used to prove our \cref{thm:bias}.

Let $\Theta$ be a probability space equipped with an appropriate $\sigma$-algebra. Consider $F: \Theta \times \R^n \to \R$ and a set $B \subset \R^n$, given by
$$
	B =\left\{{\boldsymbol \beta} : g_i({\boldsymbol \beta})=0, i \in I ; g_i({\boldsymbol \beta}) \leq 0, i \in J\right\}.
$$
where $I$ and $J$ are finite index sets.
Let ${P}$ and $Q$ be two probability measures on $\Theta$.
Let $\phi( {\boldsymbol \beta}, \alpha ) = ( {P} + \alpha(Q - {P} ) )F(\cdot, {\boldsymbol \beta}) $ and $\phi({\boldsymbol \beta}) = \phi( {\boldsymbol \beta}, 0 )$. Here $\min_B \phi$ is the main program of interest and $Q$ is a perturbation measure. The amount of perturbation is measured by $\alpha \in [0,1]$.

Let $\betast$ be the unique minimizer of $\phi({\boldsymbol \beta}) $ over $B$.
Let $J\st \subset J$ be the inequality constraints active at $\betast$, meaning $g_i(\betast)=0$ for all $i \in J\st$.
Let $\betast(\alpha)$ be the unique optimal solution to $\min _B \phi(\cd, \alpha)$.

\begin{lemma}[Theorem 1 from \citet{shapiro1990differential}]

	Define the Lagrangian function by
	$ L({\boldsymbol \beta}, \lambda, \alpha )=\phi({\boldsymbol \beta}, \alpha)+\sum_{i \in I \cup J} \lambda_i g_i({\boldsymbol \beta})$ and $L_0 ({\boldsymbol \beta}, \lambda) = L({\boldsymbol \beta}, \lambda, 0)$.
	Let $\Lambda_0$ be the set of optimal Lagrangian multipliers \footnote{$\lambda \in \Lambda_0$ iff
		$
			\nabla_x L_0\left(\beta_0, \lambda\right)=0
		$
		and $\lambda_i \geq 0, i \in J\st$ and $\lambda_i = 0$ for all $j \in J \setminus J\st$.}. Define the critical cone $C = \{ {\bm u } : {\bm u }\tp \nabla g_i (\betast) = 0, i \in I; {\bm u }\tp \nabla g_i(\betast) \leq 0$ for $i \in J \st ; {\bm u } \tp \nabla \phi(\betast) \leq 0 \}$. Define $\Lambda_0\st = \argmax_{\lambda \in \Lambda_0} \sum_{i \in I \cup J} \lambda_i g_i(\betast) $.
	Assume the following conditions.

	\label{lm:shapiro}
	\begin{enumerate}[1.]
		\item Differentiability. The functions $\phi(\cd, 0)$, $\phi(\cd, 1)$, and $g_i, i \in I \cup J$, are continuously differentiable in a neighborhood of $\betast$.
		\item Constraint Qualification. The gradients $\nabla g_i (\betast) $, $i \in I$ are linearly independent. And there exists $b \in \Rn$ such that $b \tp \nabla g_i(\betast) = 0$ for $i \in I$ and $ b \tp \nabla g_i (\betast) < 0$ for $i \in J\st$.
		\item Differentiability. The functions $\phi(\cd, 0)$, $g_i$, $i\in I \cup J\st$ are twice continuously differentiable in a neighbourhood of $\betast$.
		\item Second-order sufficient condition. Assume for all nonzero ${\bm u } \in C$, $\max_{\lambda \in \Lambda_0\st } {\bm u } \tp \nabla^2_{xx} L_0 (\betast, \lambda) {\bm u } > 0$.
	\end{enumerate}

	Define the function $b$ and the set ${\bm \Sigma}$ by
	\begin{align}
		b({\bm u })   & = \max_{\lambda \in \Lambda_0\st} {\bm u }\tp \nabla^2_{xt} L_0(\betast, \lambda) + \frac12 {\bm u } \tp \nabla^2_{xx} L_0 (\betast, \lambda) {\bm u }
		\label{eq:defb}
		\\
		{\Sigma} & = \argmin {\bm u } \tp \nabla \phi (\betast)
		\label{eq:defSigma}
		\\
		       & \mathrm{ s.t. }\quad  {\bm u } \tp \nabla g_i (\betast) = 0, i \in I; {\bm u } \tp \nabla g_i (\betast)  \leq  0, i \in J \st ; \notag
	\end{align}
	Then (i) there exists a positive $K$ such that $\| \betast(\alpha) - \betast\| \leq K \alpha $ for all positive $\alpha$ in a neighborhood of zero. (ii) If, in addition, the function $b$ has a unique minimizer ${\bm u }$ over ${\Sigma}$, then the limit
	$\lim _ {\alpha \downarrow 0} ( \betast (\alpha) - \betast) / \alpha$ exists and equals to ${\bm u }$.
\end{lemma}

\section{Proofs}

\subsection{Proof of Closed-Form Expression for Bias and \cref{thm:bias}}

\label{sec:proof:thm:bias}


Now we prove \cref{thm:bias}.

To begin with, we define directional derivatives. For a probability measure ${P}$ on the item set $\Theta$, let ${\boldsymbol \beta} ({P})$ be the unique optimal solution to the Eisenberg-Gale program: ${\boldsymbol \beta}({P}) = \min _{ {\boldsymbol \beta} \in (0, 1]^n} \int F(\theta, {\boldsymbol \beta}) \diff {P}(\theta) $.
Let ${P}$ and $Q$ be two measures.
When exists, the directional derivative is defined as
\begin{align*}
	\diff {\boldsymbol \beta} ( {P}; Q-{P}) = \lim_{t \downarrow 0} \frac{{\boldsymbol \beta} ( {P} + t(Q-{P}) ) - {\boldsymbol \beta}({P})}{t}.
\end{align*}

Step 1. Show there exist $K > 0$ and $\bar \alpha _1> 0$ such that $\| \betasta - \betast \| \leq K \alpha $ for $\alpha \in [0, \bar \alpha _ 1]$. And $\diff {\boldsymbol \beta}(s; s - s' ) = - ({\bm P}\cH_0 { \bm P}) \pinv  \bm \delta_0$.

We apply \cref{lm:shapiro} with $F(\theta, {\boldsymbol \beta}) = \max_i \vithe \betai - \sum_i b_i \log (\betai) $, $g_i ({\boldsymbol \beta}) = \betai - 1$, $B = (0,1]^n$, $\diff {P} = s \diff \theta$, $ \diff (Q - {P}) = (s-s') \diff \theta$.

By \nameref{as:constraint_qualification} in the market $\cM_0$, the set ${\Sigma}$ (defined in \cref{eq:defSigma}) becomes the plane $\{ {\bm u } : u_i = 0, i \in I_+ \}$ where we recall $I_+ = \{i:\betasti = 0\}$.
By \nameref{as:smoothness} in $\cM_0$, we know the Lagrangian multiplier of the EG program in $\cM_0$ is unique, and thus $b(\cd)$ defined in \cref{eq:defb} becomes
$b({\bm u }) = {\bm u } \tp \bm \delta_0   + \frac12 {\bm u } \tp \cH _ 0 {\bm u }$.
We conclude
the directional derivative is
$\diff {\boldsymbol \beta}( s; s-s') = \lim_{\epsilon \to 0^+} (\betast_0 - \betast_{-\epsilon})/\epsilon= - ({\bm P}\cH_0 {\bm P}) \pinv  \bm \delta_0$, where
${\bm P} = \Diag( 1 (\betasti < 1) )$ and $\bm \delta_0 = \int \nabla f (\cd, \betast) ( s- s') \diff \theta$

Step 2. For all $\alpha \geq 0$ small enough, it holds $I_\alpha = I$, where we recall $I\sa = \{ i : \betastai = 1\}$ and $I = \{ i : \betasti = 0\}$.

First we show $I \subseteq I_\alpha$, i.e., $\{ i: \betasti = 1\} \subseteq  \{ i: \betastai = 1\}$ for $\alpha$ small enough. This is saying, if a constraint $\betai \leq 1$ is strongly active in the EG program of $\cM_0$, then it is also strongly active in the EG program of the market $\cM\sa$.
This holds by Lemma 2.2 from \citet{shapiro1988sensitivity} and
the \nameref{as:constraint_qualification} condition in $\cM_0$.
Next, we show, for small $\alpha$, $I_\alpha \subseteq I$ by showing $[n] \setminus I \subseteq  [n]\setminus I_\alpha $, i.e., $\{ i: \betasti < 1\} \subseteq \{ i: \betastai < 1\}$. This holds by $\| \betasta - \betast \| \leq K \alpha$ for all $\alpha \leq \bar \alpha_1$.

Step 3. For all $\alpha \geq 0$ small enough, the directional derivative $\diff {\boldsymbol \beta}( s\sa ; s-s') = \lim_{\epsilon \to 0^+} ( \betast_{\alpha - \epsilon} - \betasta )/\epsilon$ exists and $\alpha \mapsto \diff {\boldsymbol \beta}( s\sa ; s - s')$ is continuous.

Applying \cref{lm:shapiro} with $\diff {P}  = s\sa \diff \theta$ and $ \diff (Q-{P}) = (s-s')\diff \theta $ shows that the directional derivative exists for all $\alpha$ small enough. Now we find an expression for $ \diff {\boldsymbol \beta}( s\sa ; s-s')$.
By \nameref{as:smoothness} in $\cM_0$, twice continuously differentiability of $\b \mapsto \nabla \sq \int F(\b,\t)s'\diff\t$ in a neighborhood of $\betast$, and the Lipschitzness result from step 1,
we know $\b \to \nabla \sq \int F(\t,\b)s\sa \diff \t$ is twice continuously differentiable at $\betasta$, which implies
uniqueness of Lagrangian multiplier in the market $\cM\sa$, for all $\alpha$ small enough (Lemma 2.2 from \citet{shapiro1988sensitivity}). And thus $b(\cd)$ defined in \cref{eq:defb} becomes
$b({\bm u }) = {\bm u } \tp {\boldsymbol \delta} \sa  + \frac12 {\bm u } \tp \cH\sa {\bm u }$.
Next, the \nameref{as:constraint_qualification} condition in the market $\cM_0$
implies SCS in $\cM\sa$ for all $\alpha$ small enough.
The set ${\Sigma}$ (defined in \cref{eq:defSigma}) becomes the plane $\{ {\bm u } : u_i = 0, i \in I_\alpha \}$.
The quadratic program $\min _{\Sigma} b$ has a closed-form solution: ${\bm u } = - ({\boldsymbol {\bm P}} \sa \cH\sa {\boldsymbol {\bm P}} \sa)\pinv {\boldsymbol \delta} \sa $ where ${\boldsymbol {\bm P}} \sa = \Diag( 1 (\betastai < 1) )$, ${\boldsymbol \delta} \sa = \int \nabla f(\cd,\betasta) (s-s')\diff \t$.
And we apply $I\sa = I $ for all $\alpha$ small enough from step 2, and so ${\boldsymbol {\bm P}} \sa = {\bm P}$.
We conclude
\begin{align*}
	\diff {\boldsymbol \beta}( s\sa ; s - s') = - ( {\bm P}\cH\sa {\bm P})\pinv {\boldsymbol \delta} \sa
\end{align*}
for all $\alpha$ small enough.

Next, we show that the map $\alpha \mapsto ({\bm P}\cH\sa {\bm P})\pinv {\boldsymbol \delta} \sa $ is continuous for all $\alpha$ small enough.
To see this, note $\cH\sa = \nabla\sq \int F(\theta,\betasta) s \sa \diff \theta = (1-\alpha) \nabla\sq \int F(\cd, \betasta) s \diff\theta + \alpha \nabla\sq \int F (\cd, \betasta) s' \diff \theta$ is continuous in $\alpha$. The ${\bm P}$ matrix is fixed, and so $\alpha \mapsto {\bm P} \cH\sa {\bm P}$ is continuous. Without loss of generality, suppose $\Ic = \{ i: \betasti < 0\} = [k]$ are the first $k$ buyers. Then ${\bm P} \cH \sa  {\bm P}$ creates a matrix with upper left $k$-by-$k$ block equal to the upper left $k$-by-$k$ block of $\cH \sa $, denoted $\cH_{\alpha, \Ic\Ic}$, and zeros everywhere else.
Then $({\bm P}\cH \sa {\bm P})\pinv$ is a matrix with upper left block being $(\cH_{\alpha, \Ic\Ic}) \inv$ and zeros everywhere else.
Since $\cH\sa$ is positive definite for all small $\alpha$, the submatrix $\cH_{\alpha, \Ic\Ic}$ must also be positive definite. And so $\alpha \mapsto (\cH_{\alpha,\Ic\Ic}) \inv$ is continuous, implying continuity of $\alpha \mapsto ({\bm P}\cH\sa {\bm P})\pinv $.
Finally, ${\boldsymbol \delta} \sa = \int \nabla f(\cd, \betasta) (s-s') \diff\theta$ is continuous in $\alpha$.

Step 4. The desired claim: $\betastt - \betast = o(\alpha)$.
Note $\biasformula$ defined in \cref{eq:defbias} is exactly $\diff {\boldsymbol \beta}(s\sa; s-s')$. And so
\begin{align*}
	\betastt - \betast
	 & = \betasta + \alpha \diff {\boldsymbol \beta}(s\sa; s-s') - \betast
	\\
	 & =    ( \betast - \alpha \diff {\boldsymbol \beta}(s_0; s-s') + o(\alpha) ) + \alpha \diff {\boldsymbol \beta}(s\sa; s-s') - \betast
	\\
	 & = \alpha (\diff {\boldsymbol \beta}(s\sa; s-s')  - \diff {\boldsymbol \beta}(s_0; s-s')  ) + o (\alpha ) = o(\alpha) .
\end{align*}
where the last line uses continuity in $\alpha$ of the directional derivative.
This completes the proof of \cref{thm:bias}.

\subsection{Proof of \cref{thm:betaconsistency}} \label{sec:proof:thm:betaconsistency}
Recall the estimator $\widehat {\boldsymbol \beta} = \betagama - \alpha (\widehat {\boldsymbol {\bm P}} \sa \widehat \cH\sa \widehat {\boldsymbol {\bm P}} \sa) \pinv \widehat {\boldsymbol \delta} \sa$
and the debiased surrogate is $\betastt = \betasta - \alpha ({\boldsymbol {\bm P}} \sa \cH \sa {\boldsymbol {\bm P}} \sa )\pinv {\boldsymbol \delta} \sa$.
By results from \citet{liao2023statistical} we know that \nameref{as:constraint_qualification} and \nameref{as:smoothness} in the market $\cM\sa$ and $\eta_t = o(1)$ imply $\betagama \toprob \betasta$ and $(\widehat {\boldsymbol {\bm P}} \sa \widehat \cH\sa \widehat {\boldsymbol {\bm P}} \sa) \pinv \toprob ({\boldsymbol {\bm P}} \sa \cH \sa {\boldsymbol {\bm P}} \sa )\pinv $.
Finally, $\widehat {\boldsymbol \delta} \sa \toprob {\boldsymbol \delta} \sa $ holds by the law of large numbers. We complete the proof of \cref{thm:betaconsistency}.

\subsection{Proof of \cref{thm:betaclt}} \label{sec:proof:thm:betaclt}

Recall the estimator $\widehat {\boldsymbol \beta} = \betagama - \alpha (\widehat {\boldsymbol {\bm P}} \sa \widehat \cH\sa \widehat {\boldsymbol {\bm P}} \sa) \pinv \widehat {\boldsymbol \delta} \sa$
and the debiased surrogate is $\betastt = \betasta - \alpha ({\boldsymbol {\bm P}} \sa \cH \sa {\boldsymbol {\bm P}} \sa )\pinv {\boldsymbol \delta} \sa$.
Define
\begin{align}
	 & \bm z_t = \betagama - \betasta - \alpha ({\boldsymbol {\bm P}} \sa \cH \sa {\boldsymbol {\bm P}} \sa )\pinv ( \widehat {\boldsymbol \delta} \sa - {\boldsymbol \delta} \sa )
	\\
	 & {\bm \xi} _t  = - \alpha ((\widehat {\boldsymbol {\bm P}} \sa \widehat \cH\sa \widehat {\boldsymbol {\bm P}} \sa) \pinv - ({\boldsymbol {\bm P}} \sa \cH \sa {\boldsymbol {\bm P}} \sa )\pinv   )  {\boldsymbol \delta} \sa
	\\
	 & {\bm \zeta_t} = - \alpha (\hhinv - \hinv) (\widehat {\boldsymbol \delta} \sa - {\boldsymbol \delta} \sa)
\end{align}
Then clearly $\widehat {\boldsymbol \beta} - \betastt = \bm z_t + {\bm \xi} _t + {\bm \zeta_t}$.

Step 1. We show $\sqrt t \bm z_t$ converges to a normal distribution.
By \citet{liao2023statistical}, it holds
\begin{align}
	\sqrt t (\betagama -\betasta ) = \frac{1}{\sqrt t} \sumtau \bm d_\alpha(\thetau) - \E[\bm d_\alpha] + o_p(1)
\end{align}
where $\bm d_\alpha(\theta) = - ( {\boldsymbol {\bm P}} \sa \cH\sa {\boldsymbol {\bm P}} \sa) \pinv \nabla f (\theta, \betasta)$ is the influence function.
Next, we define the likelihood ratio function $g (\theta) = (\diff  (s - s') / \diff s \sa ) (\theta) = -\frac{1}{\alpha } 1_{\Theta_\bad} (\theta) + \frac{1}{1-\alpha } 1_{\Theta _\good}(\theta) $.
Then ${\boldsymbol \delta} \sa = \int \nabla f (\cd, \betasta) (s - s')\diff \theta =\int g(\cd)\nabla  f (\cd ,\betasta)  s\sa (\cd) \diff \theta $.
Then by definition $\widehat {\boldsymbol \delta} \sa$ can be written as $\widehat {\boldsymbol \delta} \sa = P_t (g (\cd)\nabla f (\cd, \betagama))$ where $P_t h(\cd) = \frac1t \sumtau h(\thetau)$ and $\thetau$ are i.i.d.\ draws from $s\sa$.
Also let $\nu_t h = \sqrt t( P_t h - \int h s\sa \diff \theta)$.
By the decomposition,
\begin{align}
	  & \sqrt t (\widehat {\boldsymbol \delta} \sa - {\boldsymbol \delta} \sa )
	\\
	= &
	\sqrt t \bigg( P_t \big( g(\cd) \nabla f (\cd, \betagama)\big) - \int g(\cd)\nabla  f (\cd ,\betasta)  s\sa \diff \theta \bigg)
	\\ =
	  & \sqrt t \bigg( P_t \big( g(\cd) \nabla f (\cd, \betasta)\big) - \int g(\cd)\nabla  f (\cd ,\betasta)  s\sa \diff \theta \bigg)
	+ \nu_t \big(g (\cd)  ( \nabla f (\cd, \betagama) - \nabla f (\cd, \betasta) ) \big)
	\\
	= & \nu_t( g(\cd) \nabla f (\cd, \betasta)) + o_p(1),
\end{align}
where the last line follows by a stochastic equicontinuity argument as in \citet{liao2023statistical}.
Finally, $ \bm d_\alpha (\cd) - \alpha g(\cd)({\boldsymbol {\bm P}} \sa \cH \sa {\boldsymbol {\bm P}} \sa )\pinv  \nabla f (\cd, \betasta) $ is exactly the influence function $\IF_1 (\cd)$ defined in the theorem.

Step 2. Show ${\bm \xi_t} = O_p(\eta_t)$.

We need a matrix perturbation result.
\begin{lemma}[Theorem 2.2 from \citet{stewart1977perturbation}] \label{lm:matrixinverse}
	Let ${\bm A}$ and ${\bm B} = {\bm A} + {\bm E}$ be nonsingular square matrices, and $\|{\bm A} \inv \| \|{\bm E} \|< 1$.
	Here $\| \cd \|$ is the operator norm. Then
	\begin{align}
		\| {\bm B}\inv - {\bm A} \inv \| \leq  \frac{\|{\bm A} \inv \| \sq }{ 1 - \|{\bm E}\| \|{\bm A}\inv \| } \| {\bm E}\|.
	\end{align}

\end{lemma}

By \citet{liao2023statistical}, we know $\P(  \widehat {\boldsymbol {\bm P}} \sa = {\boldsymbol {\bm P}} \sa ) \to 1$.
Without loss of generality, suppose $I^c\sa  = \{ i: \betastai < 0\} = [k]$ are the first $k$ buyers.
Then $ ({\boldsymbol {\bm P}} \sa \cH \sa  {\boldsymbol {\bm P}} \sa )$ creates a matrix with upper left $k$-by-$k$ block equal to the upper left $k$-by-$k$ block of $\cH \sa $, denoted $\cH_{\alpha, \Ic\Ic}$, and zeros everywhere else.
Under the event $\{  \widehat {\boldsymbol {\bm P}} \sa = {\boldsymbol {\bm P}} \sa\}$,
the matrix $(\widehat {\boldsymbol {\bm P}} \sa \widehat \cH \sa \widehat {\boldsymbol {\bm P}} \sa)$ is one with upper left $k$-by-$k$ block being upper left $k$-by-$k$ block of $\widehat \cH \sa $, denoted $\widehat \cH_{\alpha, \Ic\Ic}$, and zeros everywhere else.
Now let ${\bm A} = \widehat \cH_{\alpha, \Ic\Ic}$ and ${\bm B} = \cH_{\alpha , \Ic\Ic}$.
It is clear that $ \| {\bm \xi_t} \|_2 \leq \alpha \|{\boldsymbol \delta} \sa\|_2 \| {\bm A} \inv - {\bm B}\inv\| = O _p (\eta_t)$ by \cref{lm:matrixinverse}.

Step 3. Show ${\bm \zeta_t} = o_p(\frac{1}{\sqrt t})$.
From previous derivation, $\widehat {\boldsymbol \delta} \sa - {\boldsymbol \delta} \sa = O_p(1/\sqrt t)$ and $ \hhinv - \hinv = O_p(\eta_t) = o_p(1)$. We conclude ${\bm \zeta_t} = O_p(1/\sqrt t) o_p(1) = o_p(1 / \sqrt t)$.

We complete the proof of \cref{thm:betaclt} part 1.

Now we prove part 2. The claim that $\cH \sa = \Diag( b_i / (\betastai)\sq )$ follows from \citet{liao2023fisher}. Next we derive the asymptotic distribution.
We only need to handle the ${\bm \xi_t}$ term defined as before. Under the event $\{ \widehat {\boldsymbol {\bm P}} \sa = {\boldsymbol {\bm P}} \sa\}$, we have ${\bm \xi_t} = - \Diag(  \alpha \delta_{\alpha ,i} 1(\betastai < 1) b_i \inv ) ( \betagama \circ \betagama - \betasta \circ \betasta)$.
By \citet{liao2023statistical} we know $\sqrt t (\betagama - \betasta) = \frac1{\sqrt t} \sumtau  (\bm d_\alpha(\thetau) - \E\sa [\bm d_\alpha]) + o_p(1)$. So by the delta method,
\begin{align*}
	\sqrt t  ( \betagama \circ \betagama - \betasta \circ \betasta) = \Diag(2 \betasta) \frac1{\sqrt t} \sumtau  (\bm d_\alpha(\thetau) - \E\sa [\bm d_\alpha]) + o_p(1).
\end{align*}
Summarizing,
\begin{align*}
	\sqrt t {\bm \xi_t} = - \Diag(  2\betastai \alpha \delta_{\alpha ,i} 1(\betastai < 1) b_i \inv ) \frac1{\sqrt t} \sumtau  (\bm d_\alpha(\thetau) - \E\sa [\bm d_\alpha]) + o_p(1)
\end{align*}
Note that for $i$ such that $\betastai = 1$, the $i$-th entry of $\bm d_\alpha$ will be zero almost surely, so the indicator can be dropped in the diagonal matrix.
Now let $\IF_2(\theta) = -\Diag(  2\alpha \betastai  \delta_{\alpha ,i} b_i \inv )  \bm d_\alpha + \IF_1$, where we recall $\IF_1$ is defined in \cref{thm:betaclt}. Then under the condition that
$\E\sa [\frac{1}{\bidgap(\betasta, \theta)}] < \infty$, we have $\sqrt t (\widehat {\boldsymbol \beta} - \betastt) \tod \cN(0, \E\sa[ (\IF_2 - \E\sa[\IF_2])  (\IF_2 - \E\sa[\IF_2]) \tp])$.

\subsection{Proof of \cref{thm:revbias}} \label{sec:proof:thm:revbias}
\begin{proof}
	Note
	\begin{align*}
		 & \bigg | \int \max_i( \vithe \betastt_i )s(\theta) \diff \theta - \int \max_i( \vithe \betast_i )s(\theta) \diff \theta \bigg |
		\\
		 & \leq
		\int | \max_i( \vithe \betastt_i )s(\theta) - \max_i( \vithe \betast_i )s(\theta) | \diff \theta
		\\
		 & \leq \max_i \sup_\theta v_i(\theta) \cdot \| \betast - \betastt \|_\infty = o(\alpha)
	\end{align*}
	as $\alpha \downarrow 0$.

	Let $P_t h(\cd) = \frac1t \sumtau h(\thetau)$ and $\thetau$ are i.i.d.\ draws from $s\sa$, and ${\bm P} h(\cd) = \int h s\sa$.
	Let $g(\cd) = \diff s / \diff s\sa = \frac{1}{1-\alpha} 1_{\Theta_\good}$.
	Also let $\nu_t h = \sqrt t( P_t h - \int h s\sa \diff \theta)$.
	We change the measure and write
	$\widetilde{\REV}\st = \int  f(\theta, \betastt)  s \diff\theta = \int f(\theta, \betastt) g(\theta) s\sa \diff\theta $.
	For the claim regarding asymptotic normality,
	note
	\begin{align*}
		 & \sqrt t (\widehat{\REV}  - \widetilde{\REV}\st )
		\\
		 & = \sqrt t [ P_t f(\cd, \widehat {\boldsymbol \beta}) g(\cd) - {P} f(\cd, \betastt)g(\cd)]
		\\
		 & =  \sqrt t  [ P_t f(\cd, \betastt)g(\cd ) - {P} f(\cd, \betastt)g(\cd ) ] +
		\sqrt t  [ {P} f (\cd, \widehat {\boldsymbol \beta}) g(\cd ) - {P} f (\cd, \betastt) g(\cd )]
		\\
		 & \quad +
		\sqrt t  (P_t - {P}) (f(\cd, \widehat {\boldsymbol \beta}) g(\cd ) - f(\cd, \betastt) g(\cd ))
		\\
		 & = \nu_t (f(\cd, \betastt) g(\cd ) ) +   \nu_t (\E\sa[ g(\cd )\nabla f(\cd, \betastt)] \tp \IF(\cd)  ) + o_p(1)
	\end{align*}
	where the last line follows by the delta method and a stochastic equicontinuity argument. We complete the proof of the theorem.
\end{proof}

\section{Ommited Maintext}
\label{sec:formal_def_finite_fppe}

\subsection{Formal Definition of Finite FPPE}
\begin{defn}[Finite FPPE, \citet{conitzer2022pacing}]
	\label{def:observed_fppe}
	Given $(b,v,\sigma, \gamma)$, the
	finite FPPE, $\oFPPE(\bm b, \bm v, \sigma , \gamma)$, is the unique tuple
	$({\boldsymbol \beta},p) \in [0,1]^n \times \R^t_+ $
	such that there exist $x_i^\tau \in [0,1]$:
(First-price) For all $\tau$, $\ptau = \max_i \betai \vitau$. For all $i$ and $\tau$, $\xitau > 0$ implies $\betai \vitau =\max_k \beta_k v_k^\tau$.
(Supply and budget feasible)  For all $i$, $ \sigma \sumtau \xitau \ptau  \leq b_i$. For all $\tau$, $\sumiton \xitau \leq 1$.
(Market clearing)  For all $\tau$, $\ptau > 0$ implies $ \sumiton \xitau = 1 $.
(No unnecessary pacing) For all $i$, $ \sigma \sumtau \xitau \ptau  < b_i$ implies $\betai = 1$.
\end{defn}
\section{Related Work}\label{ssec:relatedwork}

\textbf{Pacing equilibrium.}
Pacing is a budget management strategy where bids are scaled down by a multiplicative factor in order to smooth out and control spending.
In the first-price setting,
\citet{borgs2007dynamics} study first price auctions with budget constraints in a perturbed model, where the limit prices converge to those of an FPPE.
Building on that work, \citet{conitzer2022pacing} introduce the FPPE framework to model autobidding in repeated auctions,
and discover several properties of FPPE such as shill-proofness monotonicity properties, and a close relationship between FPPE and the quasilinear Fisher market equilibrium~\citep{chen2007note,cole2017convex}.
\citet{gao2022infinite} propose an infinite-dimensional variant of the quasilinear Fisher market, which lays the probability foundations of the current paper.
\citet{gao2021online} and \citet{liao2022dualaveraging} study online computation of the infinite-dimensional Fisher market equilibrium, and we utilize their methods to compute equilibria.
In the second-price setting, a variety of models have been proposed for modeling the outcome of budget management:
\citet{balseiro2015repeated} study budget-management in second-price auctions through a fluid mean-field approximation;
\citet{balseiro2019learning} investigate adaptive pacing strategy from buyers' perspective in a stochastic continuous setting; and
\citet{balseiro2021budget} study several budget smoothing methods including multiplicative pacing in a stochastic context.
A second-price analogue of FPPE was explored by \citet{conitzer2022multiplicative}.
Modeling statistical inference in the second-price setting is most likely harder, due to equilibrium multiplicity issues~\citep{conitzer2022multiplicative}, and hardness of even computing equilibria~\citep{chen2021complexity}.

Recently, a series of papers (e.g., \citet{aggarwal2019autobidding,deng2021towards,balseiro2021robust}) investigated an ``autobidding'' model for value-maximizing bidders, who aim to maximize their value subject to budget or return-on-spend constraints, which may include buying individual ad impressions that have negative utility. Several extensions provide results for more general settings, such as the quality of equilibria for generalized rationality models for bidders \cite{babaioff2021non}, non-uniform bidding strategies \cite{deng2023non}, and liquid welfare guarantees when budget-constrained buyers fail to reach equilibrium \cite{gaitonde2023budget}.
Statistical inference for parallel A/B testing under such models represent an interesting future direction of work.





\textbf{Interference in A/B Testing for Auction Markets.} Experimentation under interference has been extensively studied in the past decade within the context of social networks~\cite{ugander2013graph,aronow_samii2017, eckles2017design, Athey2018,LiWager2022} and two-sided marketplaces~\cite{Zigler2018, Wager2021}. Under interference, the outcome of an experiment unit, either a user in a social network or a buyer in an online marketplace, can be affected by the treatment status of other units in the networks, violating the Stable Unit Treatment Values Assumption (SUTVA) assumption commonly employed in analyzing online A/B tests~\cite{rosenbaum2007interference}. To take into account this type of interference, novel designs and analysis approaches have been proposed, e.g., cluster randomization~\cite{ugander2013graph,karrer2021network}, regression-based estimators~\cite{chin2019}, and bipartite analysis~\cite{Harshaw2021}. A/B testing in the presence of seller-side budget constraints has been discussed by \citet{Basse2016} and \citet{Liu2021} where the budget-splitting design spans the full market.
However, most of the above works do not consider the particular type of interference caused by competitive equilibrium effects.
Recent work has also addressed equilibrium effects~\citep{Wager2021,liao2023fisher,liao2024bootstrap}. Most relevant to our results, \citet{liao2023statistical} consider A/B testing via budget splitting in Fisher markets and FPPE. However, none of these works consider the effects caused by parallel A/B tests that interfere with each other.
Our method builds on top of the budget-splitting design; our approach additionally segments the full market into submarkets, where different submarkets receive different treatment. In the formal sections of our paper we assume that the submarket structure is known. In practice, the submarkets may not be known; \citet{rolnick2019randomized} and \citet{viviano2023causal} present methods to identify submarkets through clustering approaches that aim to limit the interference between clusters.


\section{Results for Revenue} \label{sec:res_rev}

\begin{algorithm}
	\caption{Debiasing procedure} \label{alg}
	\textbf{Input.}
	Budgets $[b_1,\dots, b_n]$. Values $\vitau$. First-price allocation $\xitau \in [0,1]$.
	The observed equilibrium pacing multipliers $\betagama = [\betagam_1, \dots, \betagam_n]$.
	Proportion of bad items $\alpha$. Desired confidence level $(1-c) \in (0,1)$.

	1. Estimate the Hessian matrix $\widehat \cH \sa$, e.g.
	via
	$\widehat \cH \sa = \Diag( b_i / (\betagamai)\sq)$, or using finite differences.

	2. Compute relevant quantities
	$\widehat {\boldsymbol {\bm P}} \sa$, $\widehat {\boldsymbol \delta} \sa$, $\mutau$.

	3. Compute the debiased estimator $\widehat {\boldsymbol \beta}$ and  covariance matrix $\widehat {\bm \Sigma}_1$.

	4. Output confidence region: $\widehat {\boldsymbol \beta}+ (\chi_{n, c} / \sqrt{t})\widehat {\bm \Sigma}_1^{1/2} \B$, where $\chi_{n, c}\sq$ is the $(1-c)$-quantile of a chi-square distribution with degree $n$, and $\B$ is the unit ball in $\Rn$.

\end{algorithm}

The revenue in the good market is defined as $\REVst = \int \max_i \{ \vithe \betasti \} s(\theta) \diff \theta$. We define the debiased revenue surrogate and its estimator
\begin{align}
	 & \widetilde{\REV} ^* = \int \max_i \{ \vithe \betastt_i \}s(\theta) \diff \theta
	\;,
	\label{eq:def_debias_rev}
	\\
	 & \widehat{\REV} = \frac{1}{1-\alpha}\frac1t \sumtau
	1_{\Theta_\good} (\thetau)
	\max_i \{ \vithetau \widehat {\boldsymbol \beta} _i \}
	\;.
\end{align}
In words, we approximate $\REVst$ with the revenue generated in the market where buyers bid according to the pacing profile $\betastt$, and then we estimate this via plug-in estimation on the subset of items that are good.
\begin{theorem}
	Let the conditions from \cref{thm:bias} and \cref{thm:betaclt} hold.
	\label{thm:revbias}
	The debiased revenue removes first-order bias:
	$\widetilde{\REV} ^* = \REVst + o(\alpha)$.
	Assume that $\sqrt t (\widehat {\boldsymbol \beta} - \betastt)$ has influence function $\IF$,
	i.e., $\sqrt t (\widehat {\boldsymbol \beta} - \betastt) = \frac{1}{\sqrt t}\sumtau \IF(\thetau) + o_p(1)$,
	then $\sqrt t (\widehat{\REV}  - \widetilde{\REV} ^*  ) \tod \cN(0, \sigma\sq_\REV)$ where
	$\sigma\sq_\REV =  \cov\sa (\IF_\REV)$,
	$\IF_\REV(\cd) = \frac{1}{1-\alpha} \big( 1_{\Theta_\good}(\cd) f(\cd, \betastt) + \E\sa[\nabla f(\theta, \betastt) 1_{\Theta_\good}(\theta)] \tp \IF(\cd)\big) $.
	The proof is in \cref{sec:proof:thm:revbias}.
\end{theorem}
Examples of influence functions compatible with \cref{thm:revbias} are the influence functions in \cref{thm:betaclt}.

We have two methods to construct CIs for the debiased revenue surrogate in \cref{eq:def_debiased_surr}.
The first method is based on the asymptotic normality results in \cref{thm:revbias}.
In the second method, a revenue CI is constructed based on a CI for $\betasta$.
Suppose $\CI$ is some confidence region for $\betasta$.
Then $\{ \frac{1}{1-\alpha}\frac1t \sumtau  \allowbreak   1_{\Theta_\good} (\thetau) \max_i \vithetau \betai: {\boldsymbol \beta} \in \CI \}$ serves as a CI for revenue.
Another way to utilize a \pmr CI is based on the observation that revenue is a monotone function of the pacing multiplier.  Let
$ \ubar{\betai}  = \min \{ \betai : {\boldsymbol \beta} \in \CI \}$ and $ \bar{{\boldsymbol \beta}}_i  = \max \{ \betai : {\boldsymbol \beta} \in \CI \}$.
Then another natural (potentially wider) CI for $\widetilde{\REV} ^*$ is $[ \underline{\REV}, \overline{\REV} ]$, where
\begin{equation} \label{eq: revenue_ci}
	\underline{\REV} = \frac1t \sumtau    \frac{1_{\Theta_\good} (\thetau)}{1-\alpha}  \max_i \vithetau \ubar{\betai} \; ,
	\quad
	\overline { \REV} = \frac1t \sumtau   \frac{1_{\Theta_\good} (\thetau)}{1-\alpha}  \max_i \vithetau \bar{{\boldsymbol \beta}}_i
	\;.
\end{equation}

\section{Parallel A/B Test under Interference}\label{sec:ABTest}
In this section, we use the theory from the previous sections to formulate a statistical inference theory for parallel A/B tests under interference.

Consider $K$ auction markets in which we want to run parallel A/B tests.
There are $n$ buyers, each with budget $b_i$ corresponding to one unit supply of items,
and an item set $\Theta$, partitioned into $\Theta_\good$ and $\Theta_\bad$.
Let $C_1, \dots, C_K$ be a partition of the buyers,
$\Theta_1, \dots, \Theta_K$ be a partition of the good item set $\Theta_\good$,
and $s_1,\dots s_K$ be a set of supply functions, supported on $\Theta_1,\dots, \Theta_K$ respectively.
The $k$-th submarket
consists of buyers in $C_k $, the item set $\Theta_k$ and supply $s_k$.
Let $s = \frac1K \sum_k s_k$ be the average mixture and $s'$ be a supply supported on $\Theta_\bad$.

To model treatment application we introduce the \emph{potential value functions}
$$ \hspace{-.2cm} v(0) \defeq (v_1(0,\cdot),\dots, v_n(0,\cdot)),\;v(1)\defeq (v_1(1,\cdot),\dots, v_n(1,\cdot))\;.$$
If item $\theta$ is exposed to treatment $w \in \{0,1\}$, then its value to buyer $i$ will be $v_i(w,\theta)$. 

We extend the structure on $\Theta_\bad$ and $\Theta_\good$ introduced in \cref{sec:inf_as_con} to A/B testing here.
We assume these submarkets are separated.
This is to model the ideal case where there is no interference.
A buyer $i \in C_k$ is only interested in items from the submarket he belongs to: $v_i(w, \theta) = 0$ for $\theta \in C_{k'}, k'\neq k$, and $w \in \{0,1\}$.
An item is called a bad item if it has positive values for buyers from different submarkets. Formally,
$\theta \in \Theta_\bad$ if there exist $w\in\{0,1\}$, $i \in C_k$, $j\in C_{k'}$, $k\neq k'$, so that $v_i(w, \theta) > 0$ and $v_j (w,\theta) > 0$.
A buyer $i$ from submarket $k$ has a (potentially) positive value only for items in $\Theta_k\cup \Theta_\bad$.

\textbf{Step 1. The treatment effects}. Let
$    \cM  (w, k)= \FPPE( \{ b_i \} _ { i\in C_k}, \{v_i(\cd, w) \}_{i \in C_k}, s_k, \Theta_k)
$
be the $k$-th limit submarket under treatment $w\in \{0,1\}$, and let revenue be $\REVst_k(w)$.
Then the treatment effect in the submarket $k$ is $\tau\st_k = \REVst_k(1) - \REVst_k(0)$.
Equivalently, all $K$ submarkets can be formed simultaneously,
$    \cM(w) = \FPPE \left(\frac1K \bm b  , \bm v(w), s , \Theta \right)
	.$
The scaling $1/K$ in budgets ensures that a buyer in $C_k$ has budget $b_i/K$ to bid for items in $\Theta_k$,
whose supply in $\cM(w)$ is $ \frac1K s_k$.
The $k$-th component of $\cM(w)$ corresponds to $\cM (w, k)$.

\textbf{Step 2. The experiment.}
A practical parallel A/B test scheme is as follows. 

\emph{Step 2.1. Budget splitting.}
Decide on a budget split ratio $\pi_0$, $\pi_1$ satisfying $\pi_0 + \pi_1 = 1$, and split budgets accordingly.
Choose A/B test size $t$, the total number of impressions in each submarket summed across treatments.
Then, the budget of buyer $i$ for treatment $0$ will be $\pi_0 t b_i$,
and $\pi_1 t b_i$ for treatment 1.
This is because $b_i$ is normalized to correspond to one unit of total supply.

\emph{Step 2.2. Observe markets with interference.}
Draw items/impressions from $\alpha s' + (1-\alpha )s$.
Draw $\lfloor \pi_0 t \rfloor$ items for each submarket under treatment $0$ and $\lfloor \pi_1 t \rfloor$ for treatment 1.
A total of $t_w =\lfloor \pi_w t K \rfloor$ items are drawn for the whole market under treatment $w$.
Since items come from the contaminated distribution, roughly a small fraction $\alpha$ of the $ t_w $ items cause interference.

The data observed following the described A/B test scheme can be compactly represented by two finite FPPEs (one for each treatment) (\cref{def:observed_fppe}):
$\widehat \cM \sa (w) = \oFPPE ( \pi_w t \bm  b, \bm v (w), 1, \gamma)$
,
where $\gamma \subset \Theta $ consists of $ t_w = \lfloor \pi_w    t K \rfloor$ i.i.d.\ draws from
the mixture $\alpha s' + (1-\alpha )s$.

It can be shown that
the observed market $\widehat \cM \sa (w)$ converges to
$\cM \sa (w) = \FPPE\
	\left(\frac1K \bm b  , \bm v(w), s \sa , \Theta \right)
	.$

\textbf{Step 3. Debiasing ${\boldsymbol \beta}$.}
Following previous section,
we can debias $\widehat \cM \sa (w)$ to approximate $ \cM(w) $.
Let $\widehat {\boldsymbol \beta}(w)$ be the debiased pacing multiplier.

\textbf{Step 4. Revenue estimator.} Estimate
$\widehat \REV _ k(w)$ by
\begin{align*}
	\frac{1}{\pi_w (1-\alpha) t } \sum_{\tau = 1}^{t _ w}
	1_{\Theta_k} (\thetau)
	\max_{i \in C_k}\vithetau \widehat \betai(w)
	\;.
\end{align*}
Then $\widehat \tau _k = \widehat \REV_k(1) - \widehat \REV _k(0)$ and
a variance estimate can be obtained by either
asymptotic normality, confidence region of ${\boldsymbol \beta}$, or bootstrap.

\section{Heuristic Derivation of the Debiased Surrogate} \label{sec:heuristic}

We present a heuristic argument for \cref{eq:defbias} here; for a rigorous treatment we need a perturbation theory result which is given in \cref{lm:shapiro} in the appendix.
Define the contaminated EG objective $H\sa (\b) = \int F(\b,\t)s\sa\diff\t$. Then $\cH\sa = \nabla\sq H\sa(\betasta)$. The map $(\b,\alpha)\mapsto H\sa(\b)$ has gradient $[\nabla H (\betast), \int F(\betasta,\t) (s'-s)\diff\t] \tp $ and Hessian $[\cH\sa, -{\boldsymbol \delta} \sa; -{\boldsymbol \delta} \sa\tp, 0]$ at $(\betasta, \alpha)$.
Then we have the following quadratic approximation of $H_{\alpha-\epsilon}(\b )$
\begin{align*}
	H_{\alpha-\epsilon}(\b ) \approx H\sa(\betasta) +
	\begin{pmatrix}
		\b - \betasta \\ -\epsilon
	\end{pmatrix} \tp
	\begin{pmatrix}
		\nabla H \sa (\betasta) \\
		\int F(\betasta,\t) (s'-s)\diff \theta
	\end{pmatrix}
	\\ +
	\frac12
	\begin{pmatrix}
		\b - \betasta \\ -\epsilon
	\end{pmatrix} \tp
	\begin{pmatrix}
		\cH\sa        & -{\boldsymbol \delta} \sa \\
		-{\boldsymbol \delta} \sa\tp & 0
	\end{pmatrix}
	\begin{pmatrix}
		\b - \betasta \\ -\epsilon
	\end{pmatrix}
	\;.
\end{align*}
And so
\begin{align*}
	\betast_{\alpha-\epsilon}
	= \argmin_{{\boldsymbol \beta} \in (0, 1]^n} H_{\alpha - \epsilon}(\b)
	\approx
	\argmin_{{\boldsymbol \beta} \in (0, 1]^n}  (\nabla H\sa(\betasta) +  \epsilon {\boldsymbol \delta} \sa) \tp \b \\+ \tfrac12 (\b - \betasta)\tp \cH\sa (\b - \betasta)
	\;,
\end{align*}
where we dropped terms that do not involve $\b$.

Next we need two observations. For $\epsilon$ small enough, we have $\betast_{\alpha-\epsilon, i} = 1 $ if $\betastai = 1$.
Intuitively, buyers in the market $\cM\sa$ with $\betastai=1$ have some nonzero leftover budget, which acts as a buffer to small changes in supply, and thus their pacing multipliers remain one.
Next, we also note ${\boldsymbol \beta} \mapsto \nabla H\sa (\betasta)\tp \b$ is constant on the set $\{ \b: \betai = 1 \text{ if } \betastai = 1 \}$ by complementary slackness.
Then we further simplify and obtain
\begin{align*}
	\betast_{\alpha-\epsilon} \approx \argmin_{\betai = 1 \text{ if } \betastai = 1}
	\epsilon {\boldsymbol \delta} \sa \tp \b + \tfrac12 (\b - \betasta)\tp \cH\sa (\b - \betasta)
	\;.
\end{align*}
The right-hand side is just a linearly constrained quadratic optimization problem, which admits a closed-form solution $\betasta - \epsilon({\boldsymbol {\bm P}} \sa\cH\sa {\boldsymbol {\bm P}} \sa)\pinv {\boldsymbol \delta} \sa$.

\input{appendix_full_synthetic_experiment}

%% file: appendix_full_synthetic_experiment.tex
\section{Fully synthetic experiments} \label{ssec: simulation}
\begin{figure}
	\centering
	\includegraphics[width=0.8\textwidth]{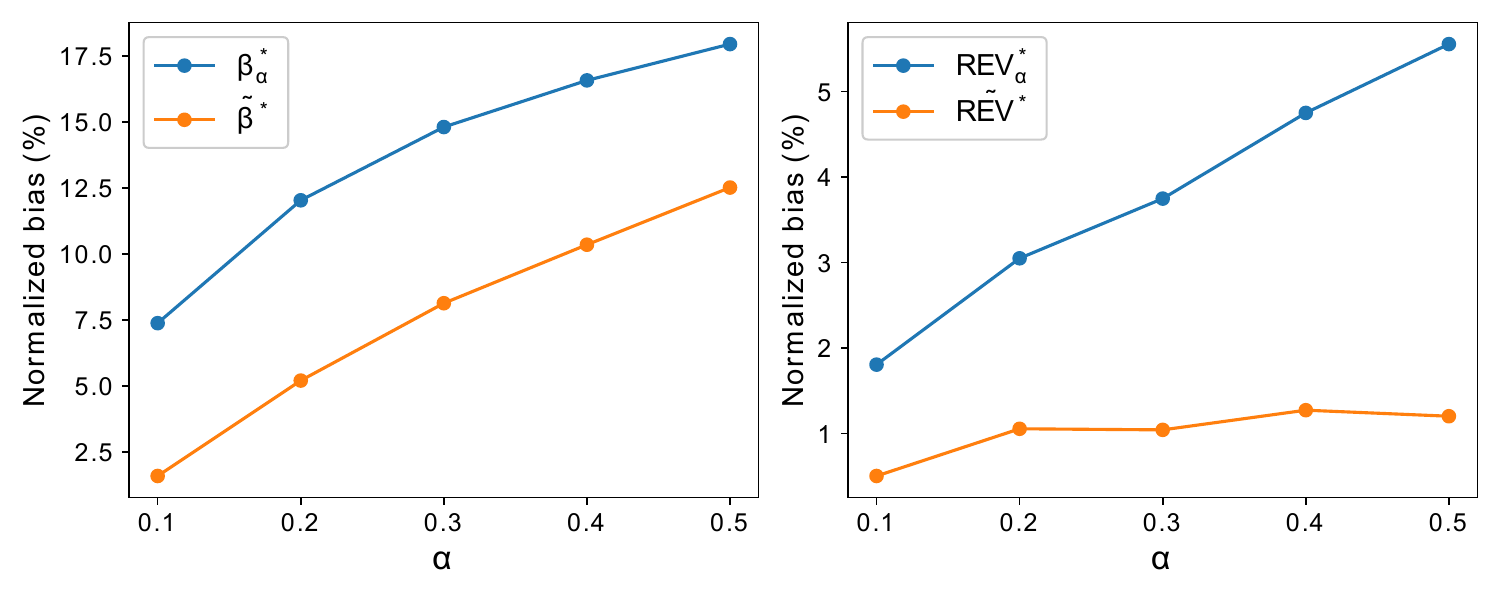}
	\caption{Normalized bias (in percent of true value) as a function of $\alpha$ in fully synthetic experiments. $\widetilde{\boldsymbol \beta}^*$ and $\widetilde\REV^*$ are the debiased surrogates for pacing multiplier and revenue in the limit market with interference $\cM\sa$. }
	\label{fig:sim_bias}
\end{figure}

In the fully synthetic experiments, we simulate 10 buyers and 1000 \emph{good} items in two submarkets, with a varying number of \emph{bad} items (up to 500) in order to study the effect of the contamination parameter $\alpha$. For each $\alpha$, we randomly sample a budget from a uniform distribution for each buyer, and compute ${\boldsymbol \beta}^*$ and $\REV^*$ from the limit pure market $\cM_0$ with a uniformly distributed value function in each submarket. \cref{fig:sim_bias} shows the normalized bias curves as a function of $\alpha$ for the fully synthetic experiment, which is similar to that from semi-synthetic experiments in~\cref{fig:semi_sim_bias}.

\begin{figure}
	\centering
	\includegraphics[width=0.5\textwidth]{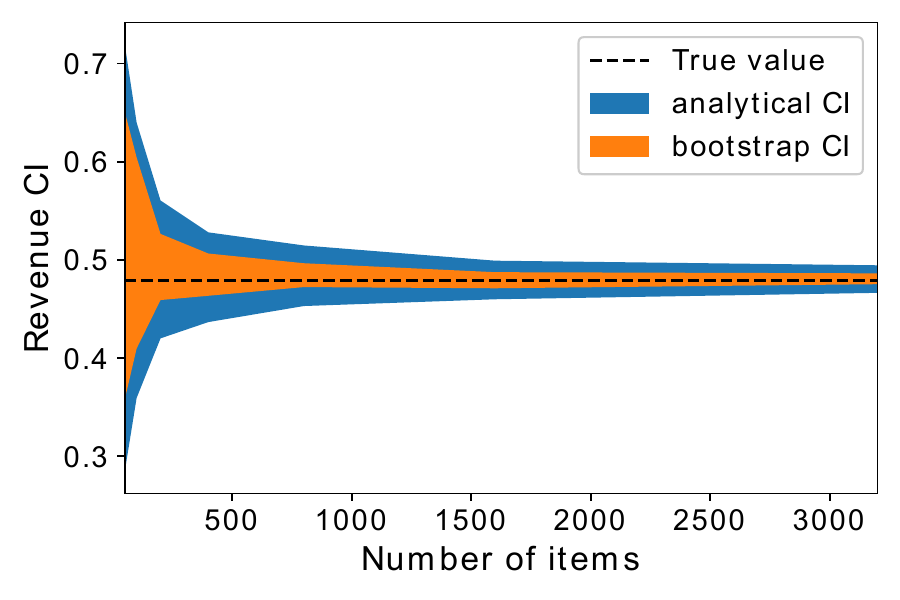}
	\caption{Revenue confidence intervals as a function of the number of items in fully synthetic experiments. The analytic CI comes from \cref{eq: revenue_ci}. The true value is the debiased surrogates for revenue in the limit market with interference $\cM\sa$. }
	\label{fig:sim_convergence}
\end{figure}

We also check the coverage of the proposed variance estimator in fully synthetic experiments. The coverage of $\widehat{\boldsymbol \beta}$ and revenue are estimated in the same way as in the semi-synthetic ones. It is shown in \cref{fig:sim_convergence} that both CIs converges to the true value from the limit market with interfrence $\cM\sa$ as the number of items goes to infinity. Then, in ~\cref{tab:full_synthetic_revenue_coverage}, we show that the bootstrap CI width is much smaller compared to that from \cref{eq: revenue_ci}. However, unlike in the semi-synthetic experiments, both CIs are over-covered, higher than the nominal level (95\%). It is expected though that the CI from \cref{eq: revenue_ci} is relatively large and hence conservative since it uses the min and max of each $\beta_i$, hence likely living outside of the confidence region ($n$-dimension ellipsoid) for $\widehat{\boldsymbol \beta}$ and leading to over-coverage. For the bootstrap CI, it requires enough bootstrap samples to fully explore the whole boundary of the confidence region of $\widehat{\boldsymbol \beta}$, as well as an accurate estimate of the asymptotic distribution of $\widehat{\boldsymbol \beta}$. In our fully synthetic simulation, our hypothesis is that the over-coverage of the bootstrap CI is mainly due to the over-estimation of the variance for $\widehat{\boldsymbol \beta}$.

\begin{table}[!ht]
	\centering\small
	\begin{tabular}{c c c c c c c c}
		\toprule
		         & \multicolumn{1}{c}{$\widehat{\boldsymbol \beta}$} & \multicolumn{4}{c}{$\widehat\REV$}                                            \\
		\cmidrule(lr){2-2} \cmidrule(lr){3-6}
		$\alpha$ & coverage                            & CI width                           & CI width (b) & coverage & coverage (b) & \\
		\midrule
		100/1100 & 0.993                               & 0.132                              & 0.043        & 1.0      & 1.0            \\
		200/1200 & 0.994                               & 0.141                              & 0.059        & 1.0      & 1.0            \\
		300/1300 & 0.980                               & 0.107                              & 0.059        & 1.0      & 1.0            \\
		400/1400 & 0.983                               & 0.098                              & 0.064        & 1.0      & 1.0            \\
		500/1500 & 0.971                               & 0.108                              & 0.069        & 1.0      & 1.0            \\
		\bottomrule
	\end{tabular}
	\caption{Coverage of $\widehat{\boldsymbol \beta}$ and revenue estimates in fully synthetic experiments. All quantities are averaged over 100 simulations for each $\alpha$. The coverage of $\widehat{\boldsymbol \beta}$ is averaged over all components of $\widehat{\boldsymbol \beta}$. For revenue estimates, the columns with "(b)" represent the quantities from the bootstrap CI, while the columns without "(b)" are for the CI from \cref{eq: revenue_ci}. The CI widths are normalized by the revenue from the limit market $\cM\sa$. $\alpha$ is expressed as the ratio of the number of bad items and the total items.}
	\label{tab:full_synthetic_revenue_coverage}
\end{table}

%% file: iclr2025_conference.bbl
\begin{thebibliography}{43}
\providecommand{\natexlab}[1]{#1}
\providecommand{\url}[1]{\texttt{#1}}
\expandafter\ifx\csname urlstyle\endcsname\relax
  \providecommand{\doi}[1]{doi: #1}\else
  \providecommand{\doi}{doi: \begingroup \urlstyle{rm}\Url}\fi

\bibitem[Aggarwal et~al.(2019)Aggarwal, Badanidiyuru, and
  Mehta]{aggarwal2019autobidding}
Gagan Aggarwal, Ashwinkumar Badanidiyuru, and Aranyak Mehta.
\newblock Autobidding with constraints.
\newblock In \emph{Web and Internet Economics: 15th International Conference,
  WINE 2019, New York, NY, USA, December 10--12, 2019, Proceedings 15}, pp.\
  17--30. Springer, 2019.

\bibitem[Aronow \& Samii(2017)Aronow and Samii]{aronow_samii2017}
Peter~M. Aronow and Cyrus Samii.
\newblock {Estimating average causal effects under general interference, with
  application to a social network experiment}.
\newblock \emph{The Annals of Applied Statistics}, 11\penalty0 (4):\penalty0
  1912 -- 1947, 2017.
\newblock \doi{10.1214/16-AOAS1005}.
\newblock URL \url{https://doi.org/10.1214/16-AOAS1005}.

\bibitem[Athey et~al.(2018)Athey, Eckles, and Imbens]{Athey2018}
Susan Athey, Dean Eckles, and Guido~W. Imbens.
\newblock Exact p-values for network interference.
\newblock \emph{Journal of the American Statistical Association}, 113\penalty0
  (521):\penalty0 230--240, 2018.
\newblock \doi{10.1080/01621459.2016.1241178}.
\newblock URL \url{https://doi.org/10.1080/01621459.2016.1241178}.

\bibitem[Babaioff et~al.(2021)Babaioff, Cole, Hartline, Immorlica, and
  Lucier]{babaioff2021non}
Moshe Babaioff, Richard Cole, Jason Hartline, Nicole Immorlica, and Brendan
  Lucier.
\newblock Non-quasi-linear agents in quasi-linear mechanisms.
\newblock In \emph{12th Innovations in Theoretical Computer Science Conference
  (ITCS 2021)}, volume 185, pp.\ ~84, 2021.

\bibitem[Balseiro et~al.(2017)Balseiro, Kim, Mahdian, and
  Mirrokni]{balseiro2017budget}
Santiago Balseiro, Anthony Kim, Mohammad Mahdian, and Vahab Mirrokni.
\newblock Budget management strategies in repeated auctions.
\newblock In \emph{Proceedings of the 26th International Conference on World
  Wide Web}, pp.\  15--23, 2017.

\bibitem[Balseiro et~al.(2021{\natexlab{a}})Balseiro, Deng, Mao, Mirrokni, and
  Zuo]{balseiro2021robust}
Santiago Balseiro, Yuan Deng, Jieming Mao, Vahab Mirrokni, and Song Zuo.
\newblock Robust auction design in the auto-bidding world.
\newblock \emph{Advances in Neural Information Processing Systems},
  34:\penalty0 17777--17788, 2021{\natexlab{a}}.

\bibitem[Balseiro et~al.(2021{\natexlab{b}})Balseiro, Kim, Mahdian, and
  Mirrokni]{balseiro2021budget}
Santiago Balseiro, Anthony Kim, Mohammad Mahdian, and Vahab Mirrokni.
\newblock Budget-management strategies in repeated auctions.
\newblock \emph{Operations research}, pp.\  859--876, 2021{\natexlab{b}}.

\bibitem[Balseiro \& Gur(2019)Balseiro and Gur]{balseiro2019learning}
Santiago~R Balseiro and Yonatan Gur.
\newblock Learning in repeated auctions with budgets: Regret minimization and
  equilibrium.
\newblock \emph{Management Science}, 65\penalty0 (9):\penalty0 3952--3968,
  2019.

\bibitem[Balseiro et~al.(2015)Balseiro, Besbes, and
  Weintraub]{balseiro2015repeated}
Santiago~R Balseiro, Omar Besbes, and Gabriel~Y Weintraub.
\newblock Repeated auctions with budgets in ad exchanges: Approximations and
  design.
\newblock \emph{Management Science}, 61\penalty0 (4):\penalty0 864--884, 2015.

\bibitem[Basse et~al.(2016)Basse, Azari~Soufiani, and Lambert]{Basse2016}
Guillaume~W. Basse, Hossein Azari~Soufiani, and Diane Lambert.
\newblock Randomization and the pernicious effects of limited budgets on
  auction experiments.
\newblock In Arthur Gretton and Christian~C. Robert (eds.), \emph{Proceedings
  of the 19th International Conference on Artificial Intelligence and
  Statistics}, volume~51 of \emph{Proceedings of Machine Learning Research},
  pp.\  1412--1420, Cadiz, Spain, 09--11 May 2016. PMLR.
\newblock URL \url{https://proceedings.mlr.press/v51/basse16b.html}.

\bibitem[Borgs et~al.(2007)Borgs, Chayes, Immorlica, Jain, Etesami, and
  Mahdian]{borgs2007dynamics}
Christian Borgs, Jennifer Chayes, Nicole Immorlica, Kamal Jain, Omid Etesami,
  and Mohammad Mahdian.
\newblock Dynamics of bid optimization in online advertisement auctions.
\newblock In \emph{Proceedings of the 16th international conference on World
  Wide Web}, pp.\  531--540, 2007.

\bibitem[Chen et~al.(2007)Chen, Ye, and Zhang]{chen2007note}
Lihua Chen, Yinyu Ye, and Jiawei Zhang.
\newblock A note on equilibrium pricing as convex optimization.
\newblock In \emph{International Workshop on Web and Internet Economics}, pp.\
  7--16. Springer, 2007.

\bibitem[Chen et~al.(2023)Chen, Kroer, and Kumar]{chen2021complexity}
Xi~Chen, Christian Kroer, and Rachitesh Kumar.
\newblock The complexity of pacing for second-price auctions.
\newblock \emph{Mathematics of Operations Research}, 2023.

\bibitem[Chin(2019)]{chin2019}
Alex Chin.
\newblock Regression adjustments for estimating the global treatment effect in
  experiments with interference.
\newblock \emph{Journal of Causal Inference}, 7, 05 2019.
\newblock \doi{10.1515/jci-2018-0026}.

\bibitem[Cole et~al.(2017)Cole, Devanur, Gkatzelis, Jain, Mai, Vazirani, and
  Yazdanbod]{cole2017convex}
Richard Cole, Nikhil~R Devanur, Vasilis Gkatzelis, Kamal Jain, Tung Mai,
  Vijay~V Vazirani, and Sadra Yazdanbod.
\newblock Convex program duality, fisher markets, and {Nash} social welfare.
\newblock In \emph{18th ACM Conference on Economics and Computation, EC 2017}.
  Association for Computing Machinery, Inc, 2017.

\bibitem[Conitzer et~al.(2022{\natexlab{a}})Conitzer, Kroer, Panigrahi,
  Schrijvers, Stier-Moses, Sodomka, and Wilkens]{conitzer2022pacing}
Vincent Conitzer, Christian Kroer, Debmalya Panigrahi, Okke Schrijvers,
  Nicolas~E Stier-Moses, Eric Sodomka, and Christopher~A Wilkens.
\newblock Pacing equilibrium in first price auction markets.
\newblock \emph{Management Science}, 2022{\natexlab{a}}.

\bibitem[Conitzer et~al.(2022{\natexlab{b}})Conitzer, Kroer, Sodomka, and
  Stier-Moses]{conitzer2022multiplicative}
Vincent Conitzer, Christian Kroer, Eric Sodomka, and Nicolas~E Stier-Moses.
\newblock Multiplicative pacing equilibria in auction markets.
\newblock \emph{Operations Research}, 70\penalty0 (2):\penalty0 963--989,
  2022{\natexlab{b}}.

\bibitem[Deng et~al.(2021)Deng, Mao, Mirrokni, and Zuo]{deng2021towards}
Yuan Deng, Jieming Mao, Vahab Mirrokni, and Song Zuo.
\newblock Towards efficient auctions in an auto-bidding world.
\newblock In \emph{Proceedings of the Web Conference 2021}, pp.\  3965--3973,
  2021.

\bibitem[Deng et~al.(2023)Deng, Mao, Mirrokni, Teng, and Zuo]{deng2023non}
Yuan Deng, Jieming Mao, Vahab Mirrokni, Yifeng Teng, and Song Zuo.
\newblock Non-uniform bid-scaling and equilibria for different auctions: An
  empirical study.
\newblock \emph{arXiv preprint arXiv:2311.10679}, 2023.

\bibitem[Eckles et~al.(2017)Eckles, Karrer, and Ugander]{eckles2017design}
Dean Eckles, Brian Karrer, and Johan Ugander.
\newblock Design and analysis of experiments in networks: Reducing bias from
  interference.
\newblock \emph{Journal of Causal Inference}, 5\penalty0 (1):\penalty0
  20150021, 2017.
\newblock \doi{doi:10.1515/jci-2015-0021}.
\newblock URL \url{https://doi.org/10.1515/jci-2015-0021}.

\bibitem[Gaitonde et~al.(2023)Gaitonde, Li, Light, Lucier, and
  Slivkins]{gaitonde2023budget}
Jason Gaitonde, Yingkai Li, Bar Light, Brendan Lucier, and Aleksandrs Slivkins.
\newblock Budget pacing in repeated auctions: Regret and efficiency without
  convergence.
\newblock In \emph{14th Innovations in Theoretical Computer Science Conference
  (ITCS 2023)}, volume 251, pp.\ ~52, 2023.

\bibitem[Gao \& Kroer(2022)Gao and Kroer]{gao2022infinite}
Yuan Gao and Christian Kroer.
\newblock Infinite-dimensional fisher markets and tractable fair division.
\newblock \emph{Operations Research (forthcoming)}, 2022.

\bibitem[Gao et~al.(2021)Gao, Kroer, and Peysakhovich]{gao2021online}
Yuan Gao, Christian Kroer, and Alex Peysakhovich.
\newblock Online market equilibrium with application to fair division.
\newblock \emph{arXiv preprint arXiv:2103.12936}, 2021.

\bibitem[Harshaw et~al.(2021)Harshaw, S\"avje, Eisenstat, Mirrorkni, and
  Pouget-Abadie]{Harshaw2021}
Christopher Harshaw, Fredrik S\"avje, David Eisenstat, Vahab Mirrorkni, and
  Jean Pouget-Abadie.
\newblock Design and analysis of bipartite experiments under a linear
  exposure-response model.
\newblock \emph{arXiv}, 2103.06392v3, 12 2021.

\bibitem[Hong et~al.(2015)Hong, Mahajan, and Nekipelov]{hong2015extremum}
Han Hong, Aprajit Mahajan, and Denis Nekipelov.
\newblock Extremum estimation and numerical derivatives.
\newblock \emph{Journal of Econometrics}, 188\penalty0 (1):\penalty0 250--263,
  2015.

\bibitem[Karrer et~al.(2021)Karrer, Shi, Bhole, Goldman, Palmer, Gelman,
  Konutgan, and Sun]{karrer2021network}
Brian Karrer, Liang Shi, Monica Bhole, Matt Goldman, Tyrone Palmer, Charlie
  Gelman, Mikael Konutgan, and Feng Sun.
\newblock Network experimentation at scale.
\newblock In \emph{Proceedings of the 27th acm sigkdd conference on knowledge
  discovery \& data mining}, pp.\  3106--3116, 2021.

\bibitem[Li \& Wager(2022)Li and Wager]{LiWager2022}
Shuangning Li and Stefan Wager.
\newblock {Random graph asymptotics for treatment effect estimation under
  network interference}.
\newblock \emph{The Annals of Statistics}, 50\penalty0 (4):\penalty0 2334 --
  2358, 2022.
\newblock \doi{10.1214/22-AOS2191}.
\newblock URL \url{https://doi.org/10.1214/22-AOS2191}.

\bibitem[Liao \& Kroer(2023)Liao and Kroer]{liao2023statistical}
Luofeng Liao and Christian Kroer.
\newblock Statistical inference and a/b testing for first-price pacing
  equilibria.
\newblock \emph{International Conference on Machine Learning}, 2023.

\bibitem[Liao \& Kroer(2024)Liao and Kroer]{liao2024bootstrap}
Luofeng Liao and Christian Kroer.
\newblock Bootstrapping fisher market equilibrium and first-price pacing
  equilibrium.
\newblock In \emph{International Conference on Machine Learning}, 2024.

\bibitem[Liao et~al.(2022)Liao, Gao, and Kroer]{liao2022dualaveraging}
Luofeng Liao, Yuan Gao, and Christian Kroer.
\newblock Nonstationary dual averaging and online fair allocation.
\newblock \emph{Advances in Neural Information Processing Systems},
  35:\penalty0 37159--37172, 2022.

\bibitem[Liao et~al.(2023)Liao, Gao, and Kroer]{liao2023fisher}
Luofeng Liao, Yuan Gao, and Christian Kroer.
\newblock Statistical inference for fisher market equilibrium.
\newblock In \emph{The Eleventh International Conference on Learning
  Representations}, 2023.
\newblock URL \url{https://openreview.net/forum?id=KemSBwOYJC}.

\bibitem[Liu et~al.(2021)Liu, Mao, and Kang]{Liu2021}
Min Liu, Jialiang Mao, and Kang Kang.
\newblock Trustworthy and powerful online marketplace experimentation with
  budget-split design.
\newblock In \emph{Proceedings of the 27th ACM SIGKDD Conference on Knowledge
  Discovery \& Data Mining}, KDD '21, pp.\  3319–3329, New York, NY, USA,
  2021. Association for Computing Machinery.
\newblock ISBN 9781450383325.
\newblock \doi{10.1145/3447548.3467193}.
\newblock URL \url{https://doi.org/10.1145/3447548.3467193}.

\bibitem[Rolnick et~al.(2019)Rolnick, Aydin, Pouget-Abadie, Kamali, Mirrokni,
  and Najmi]{rolnick2019randomized}
David Rolnick, Kevin Aydin, Jean Pouget-Abadie, Shahab Kamali, Vahab Mirrokni,
  and Amir Najmi.
\newblock Randomized experimental design via geographic clustering.
\newblock In \emph{Proceedings of the 25th ACM SIGKDD International Conference
  on Knowledge Discovery \& Data Mining}, pp.\  2745--2753, 2019.

\bibitem[Rosenbaum(2007)]{rosenbaum2007interference}
Paul~R Rosenbaum.
\newblock Interference between units in randomized experiments.
\newblock \emph{Journal of the American Statistical Association}, 102\penalty0
  (477):\penalty0 191--200, 2007.
\newblock \doi{10.1198/016214506000001112}.
\newblock URL \url{https://doi.org/10.1198/016214506000001112}.

\bibitem[Shapiro(1988)]{shapiro1988sensitivity}
Alexander Shapiro.
\newblock Sensitivity analysis of nonlinear programs and differentiability
  properties of metric projections.
\newblock \emph{SIAM Journal on Control and Optimization}, 26\penalty0
  (3):\penalty0 628--645, 1988.

\bibitem[Shapiro(1990)]{shapiro1990differential}
Alexander Shapiro.
\newblock On differential stability in stochastic programming.
\newblock \emph{Mathematical Programming}, 47:\penalty0 107--116, 1990.

\bibitem[Shapiro et~al.(2021)Shapiro, Dentcheva, and
  Ruszczynski]{shapiro2021lectures}
Alexander Shapiro, Darinka Dentcheva, and Andrzej Ruszczynski.
\newblock \emph{Lectures on stochastic programming: modeling and theory}.
\newblock SIAM, 2021.

\bibitem[Stewart(1977)]{stewart1977perturbation}
Gilbert~W Stewart.
\newblock On the perturbation of pseudo-inverses, projections and linear least
  squares problems.
\newblock \emph{SIAM review}, 19\penalty0 (4):\penalty0 634--662, 1977.

\bibitem[Ugander et~al.(2013)Ugander, Karrer, Backstrom, and
  Kleinberg]{ugander2013graph}
Johan Ugander, Brian Karrer, Lars Backstrom, and Jon Kleinberg.
\newblock Graph cluster randomization: Network exposure to multiple universes.
\newblock In \emph{Proceedings of the 19th ACM SIGKDD International Conference
  on Knowledge Discovery and Data Mining}, KDD '13, pp.\  329–337, New York,
  NY, USA, 2013. Association for Computing Machinery.
\newblock ISBN 9781450321747.
\newblock \doi{10.1145/2487575.2487695}.
\newblock URL \url{https://doi.org/10.1145/2487575.2487695}.

\bibitem[van~der Vaart(2000)]{van2000asymptotic}
Aad~W van~der Vaart.
\newblock \emph{Asymptotic statistics}, volume~3.
\newblock Cambridge university press, 2000.

\bibitem[Viviano et~al.(2023)Viviano, Lei, Imbens, Karrer, Schrijvers, and
  Shi]{viviano2023causal}
Davide Viviano, Lihua Lei, Guido Imbens, Brian Karrer, Okke Schrijvers, and
  Liang Shi.
\newblock Causal clustering: design of cluster experiments under network
  interference.
\newblock \emph{arXiv preprint arXiv:2310.14983}, 2023.

\bibitem[Wager \& Xu(2021)Wager and Xu]{Wager2021}
Stefan Wager and Kuang Xu.
\newblock Experimenting in equilibrium.
\newblock \emph{Management Science}, 67\penalty0 (11):\penalty0 6694–6715,
  2021.
\newblock \doi{10.1287/mnsc.2020.3844}.
\newblock URL \url{https://doi.org/10.1287/mnsc.2020.3844}.

\bibitem[Zigler \& Papadogeorgou(2018)Zigler and Papadogeorgou]{Zigler2018}
Corwin~M. Zigler and Georgia Papadogeorgou.
\newblock Bipartite causal inference with interference.
\newblock \emph{arXiv}, 1807.08660v1, 07 2018.

\end{thebibliography}
